\documentclass{amsart}

\input xy
\xyoption{all}

\usepackage{stmaryrd}
\usepackage{setspace}
\usepackage{geometry}
\usepackage{amsmath}
\usepackage{amssymb}
\usepackage{amsthm}  % see geometry.pdf on how to lay out the page. There's lots.
\usepackage{cite}
\usepackage[final]{showkeys}
\usepackage{hyperref}
\newtheorem{same}{This should never appear}[section]
\newtheorem{defin}[same]{Definition}
\newtheorem{claim}[same]{Claim}
\newtheorem{remark}[same]{Remark}
\newtheorem{theorem}[same]{Theorem}
\newtheorem{example}[same]{Example}
\newtheorem{lemma}[same]{Lemma}
\newtheorem{fact}[same]{Fact}

\newtheorem{cor}[same]{Corollary}
\newtheorem{prop}[same]{Proposition}

\newtheorem{observation}[same]{Observation}
\newbox\noforkbox \newdimen\forklinewidth
\setbox0\hbox{$\textstyle\smile$}
\setbox1\hbox to \wd0{\hfil\vrule width \forklinewidth depth-2pt
 height 10pt \hfil}
\setbox\noforkbox\hbox{\lower 2pt\box1\lower 2pt\box0\relax}
\def\unionstick{\mathop{\copy\noforkbox}\limits}

\def\nonfork_#1{\unionstick_{\textstyle #1}}

\setbox0\hbox{$\textstyle\smile$}
\setbox1\hbox to \wd0{\hfil{\sl /\/}\hfil}
\setbox2\hbox to \wd0{\hfil\vrule height 10pt depth -2pt width
              \forklinewidth\hfil}
\newbox\doesforkbox
\setbox\doesforkbox\hbox{\lower 2pt\box1 \lower 2pt\box2\lower2pt\box0\relax}
\def\nunionstick{\mathop{\copy\doesforkbox}\limits}

\def\fork_#1{\nunionstick_{\textstyle #1}}

\newcommand{\dnfg}{{}^{g}\unionstick}
\newcommand{\dnft}{{}^{t}\unionstick}

\newcommand{\bc}{\bold{c}}

\newcommand{\eff}{\mathcal{F}}

\newcommand{\bn}{\bold{n}}
\newcommand{\bx}{\bold{x}}
\newcommand{\by}{\bold{y}}
\newcommand{\bz}{\bold{z}}

\newcommand{\rest}{\upharpoonright}

\newcommand{\lcm}{\te{lcm}}

\newcommand{\bm}{\mathbf{m}}
\newcommand{\br}{\bold{r}}
\newcommand{\divides}{\mid}
\newcommand{\bv}{\bold{v}}
\newcommand{\Mod}{\text{Mod}}

\newcommand{\comment}[1]{}
\newcommand{\footnotei}[1]{}

\newcommand{\seq}[1]{\langle #1 \rangle}
\newcommand{\te}[1]{\textrm{#1}}

\title{The $\Gamma$-ultraproduct and Averageable Classes}
\author{Will Boney}
\email{wboney@math.harvard.edu}
\address{Mathematics Department\\Harvard University\\Cambridge, MA, USA}

\date{\today
This work was begun while the author was working towards his PhD under Rami Grossberg and he is grateful for his guidance and support.  Part of this material is based upon work done while the author was supported by the National Science Foundation under Grant No. DMS-1402191.} % delete this line to display the current date

\begin{document}

\maketitle

%\tableofcontents

\begin{abstract}
This paper introduces the $\Gamma$-ultraproduct, which is designed to take a collection of structures omitting some fixed set of unary types $\Gamma$ and average them into a structure that also omits those types.  The motivation comes from the Banach space ultraproduct, and generalizes other existing constructions such as the torsion submodule.  Motivated by examples and counterexamples, we explore conditions on classes that make the $\Gamma$-ultraproduct well-behaved and apply results from the existing literature on classification for nonelementary classes.  We use torsion modules over PIDs as an extended example.
\end{abstract}

\section{Introduction}

Ultraproducts are an invaluable tool in first-order model theory.  The ability to create a new structure that is the ``average'' of some collection of structures has far-reaching implications, perhaps most importantly the compactness theorem.  When trying to adapt first-order results (such as various results of classification theory) to nonelementary contexts, the lack of compactness is a major stumbling block.  One strategy is to use set-theoretic hypotheses to allow very complete ultrafilters (see \cite{makkaishelah} \cite{kolmanshelah} \cite{tamelc}).  Another approach is to assume that the class satisfies some fragment of compactness; an example of this is the property tameness, which was introduced by Grossberg and VanDieren \cite{tamenessone} and has seen a large amount of activity in recent years.

This paper takes a different approach.  Rather than appealing to the uniform construction of the ultraproduct, we fix a collection $\Gamma$ of types to be omitted in advance and then build the $\Gamma$-ultraproduct with the express purpose of creating an average that omits those types.  More precisely, fix the following:
\begin{itemize}
	\item a language $L$;
	\item a collection of unary $L$-types $\Gamma$; and
	\item a collection of $L$-structures $\{M_i :i \in I\}$ that each omit\footnote{Of course, it is only necessary that this omission happen on a measure one set.} every type in $\Gamma$.
\end{itemize}
We sometimes refer to this collection as \emph{the data}.

The main definition of this paper is the $\Gamma$-ultraproduct of this data by some ultrafilter $U$, which is denoted $\prod^\Gamma M_i/U$.
\begin{defin} \label{nkoudef}
Fix an ultrafilter $U$ on $I$.

	 \begin{eqnarray*}
	\prod^{\Gamma}_{i \in I} M_i := \{ f \in \prod M_i &:& \te{ there is some } X_f \in U \te{ such that, for each } p \in \Gamma,\te{ there }\\ 
	& & \te{ is a } \phi^p_f(x) \in p \te{ such that } M_i \models  \neg \phi^p_f\left(f(i)\right) \te{ for each } i \in X_f\}
	\end{eqnarray*}
\begin{itemize}

	\item Form $\prod^\Gamma M_i /U$ by giving it universe $\prod^\Gamma M_i/U := \{ [f]_U: f \in \prod^\Gamma M_i\}$ and inheriting the functions and relations from the full ultraproduct $\prod M_i/U$.
	
	\item If $\Gamma = \{p\}$, then we write $\prod^p M_i/U$.
	
\end{itemize}
\end{defin}
Given $[f]_U \in \prod^\Gamma M_i/U$, we call a choice function $p \in \Gamma \mapsto \phi^p_f$ as in the definition a \emph{witness} for $[f]_U$'s inclusion in the $\Gamma$-ultraproduct.  We typically denote witnesses and choice functions by $\mathcal{C}$.  Note that there are often many witnesses for a single element.

Unlike the normal ultraproduct, there is no reason to suspect that the $\Gamma$-ultraproduct is always well-behaved or is even an $L$-structure; these issues and examples of where things go wrong are explored in Section \ref{propgamup-sec}.  Note that the assumption that the types of $\Gamma$ are unary is crucial for this definition.  This allows the inclusion criteria for $\prod^\Gamma_{i \in I} M_i$ to be local, i. e., only depend on the function on consideration.  If the types were \emph{not} unary (such as those expressing a group is locally finite), then determining wether a choice function $f$ were to be included would require a witness that involves the other included choice functions in some way.  There seems to be no uniform way to generalize the above definition in this case.  Indeed, simply coding tuples as single elements with projection functions does not avoid this necessity: the resulting $\Gamma$-ultraproduct might omit the coded types, but fail to satisfy the sentences stating that finite tuples are coded as elements.

Section \ref{aecsection} introduces \emph{averageable classes} (roughly nonelementary classes where the appropriate $\Gamma$ is well-behaved) and applies some results from the classification theory of AECs; Theorem \ref{divline-thm} here gives a dividing line in the number of models for averageable classes.  Section \ref{example-sec} gives several examples of these classes, including dense linearly ordered groups with a cofinal $\mathbb{Z}$-chain\comment{ and differentially closed fields of characteristic 0 where every element is differentially algebraic over the constants}.  Section \ref{tormod-sec} develops the example of torsion modules over a PID, including the appropriate \L o\'{s}' Theorem and some stability theory. 

This construction can be seen as a generalization of two well-known constructions: ultraproducts of multi-sorted structures and Banach space ultraproducts.  Subsection \ref{sort-ex} shows how to view the ultraproduct of multi-sorted structures as the appropriate $\Gamma$-ultraproduct.  Moreover, if $\Gamma$ is finite (as it is in most of our examples, with Banach spaces and Archimedean fields being the only non-examples in this paper), there is a single type $p_\Gamma$ such that omitting $p_\Gamma$ is equivalent to omitting all of $\Gamma$.  Then, we could attempt to impose a sorted structure on a model omitting $\Gamma$ by which formula of $p_\Gamma$ it omits, and attempt to translate the language and syntax to a sorted one.  This would be an alternate presentation of these results: being able to sort the language corresponds to $\Gamma$-closed (Definition \ref{gamclose-def}) and being able to sort the formulas corresponds to $\Gamma$-nice (Definition \ref{gamnice-def}).  We chose the current presentation in part because some choice of equally valid presentations must be made, but also to accommodate cases of omitting infinitely many types and to avoid the unnaturality discussed below. We discuss a third possible presentation in Section \ref{diffapp-sec}.

In the standard Banach space ultraproduct, the elements of $\prod \mathcal{B}_i/U$ are sequences of bounded norm that are modded out by the equivalence relation
$$(x_i) \sim_U (y_i) \iff \lim_U \| x_i - y_i\|^{\mathcal{B}_i} = 0$$
In the standard model-theoretic ultraproduct, the elements of $\prod M_i/U$ are sequences that are modded out by the equivalence relation
$$(x_i) \sim_U (y_i) \iff \{ i : x_i = y_i\} \in U$$
Both constructions contain a step that ignores $U$-small differences; this is the equivalence relations.  However, the Banach space ultraprpoduct contains an extra step the excludes unbounded sequences.  In model theoretic language, this amounts to excluding sequences that would realize the type $\{ \|x\| > n : n < \omega\}$.  The model theoretic ultraproduct has no similar step.  We do this to arrive at $\prod^\Gamma M_i / U$.  Example \ref{bs-ex} goes into greater detail about the application of the $\Gamma$-ultraproduct to Banach spaces.  Indeed, this example was the original motivation for the $\Gamma$-ultraproduct as an attempt to generalize the extra step of throwing away unbounded elements to more general situations.  Note that, in continuous first order logic (see \cite{fourguys}), attention is restricted to uniformly bounded metric spaces and, thus, avoid the extra step.  Ben Yaacov \cite{byunbounded} has explored continuous logic in unbounded metric spaces.  The key there is to restrict the logic to only allow quantifiers that specify that the type is omitted, that is, quantifiers that turn $\phi(x, \by)$ into $\exists x \left(\|x\|< n \wedge \phi(x, \by)\right)$ for some $n <\omega$; see Observation \ref{doublestar} for a discussion of that technique here..

Many of the proofs of this paper (particularly the basic exploration of the $\Gamma$-ultraproduct in Section \ref{propgamup-sec}) are straightforward (especially in light of the above comparisons).  However, there seems to be no place in the literature that discusses these constructions in this generality or applies them to achieve compactness-like results in nonelementary classes.  In particular, the results of Section \ref{tormod-sec} on the compactness and classification theory of torsion modules over PIDs is new.  While these results could have been obtained by ``sorting" the structures and applying the ultraproduct of sorted structures\footnote{Or, in the case $R = \mathbb{Z}$, imposing a metric on torsion abelian groups by setting $d(g, h) = \log o(g-h)$.}, this class is always considered as a nonelementary (single-sorted) class.  Moreover, the translation to a sorted class would be very unnatural: the single addition function $+$ would replaced by  a collection of addition functions $\{+{r,r'} \mid r, r'\in R\}$ (and similarly for other functions).  Moreover, formulas like $``\exists z(x+z = y)"$ would not survive the sorting translation and one would be forced to specify an annihilator of $x$, $y$, and $z$ to use such a formula.  Thus, we prefer to work with torsion modules as a ``sortable" class, rather than one that is actually sorted.

\section{Properties of $\prod^\Gamma M_i/U$} \label{propgamup-sec}

Our main goal will be analyzing compactness in classes of the form $\left(EC(T, \Gamma), \prec\right)$ or $\left(EC(T, \Gamma), \subset\right)$ via the $\Gamma$-ultraproduct (recall the definition of the $\Gamma$-ultraproduct from Definition \ref{nkoudef} and that $EC(T, \Gamma)$ is the class of all models of $T$ that omit each type in $\Gamma$).  However, in this section we analyze this construction in more generality; we specialize back to these classes in Section \ref{avclasses}.

This definition and the discussion below work in a great deal of generality, in particular allowing many unary types of different sizes.  In concrete cases, $\Gamma$ often consists of a single countable type and the reader can simplify to this case with little loss.

The analysis of $\prod^\Gamma M_i/U$ breaks along two main questions:
\begin{itemize}
	\item Is $\prod^\Gamma M_i/U$ a structure, specifically a substructure of $\prod M_i/U$?
	\item Is $\prod^\Gamma M_i/U$ an elementary substructure of $\prod M_i/U$?
\end{itemize}

We analyze each of this separately, although we first provide examples that the answer to each question can be no.

\begin{example} \ \label{badex}
\begin{enumerate}
	\item Set $M = (\omega, +, \mid, 2)$, $I = \omega$, $p(x) = \{ (2^k \mid x) \wedge (x \neq 0) : k < \omega\}$, where the $`\mid'$ is the symbol for `divides.'  Then $[n \mapsto 1]_U, [n\mapsto 2^n - 1]_U \in \prod^p M/U$, but 
$$[n \mapsto 1]_U + [n\mapsto 2^n - 1]_U  = [n \mapsto 2^n]_U \not\in\prod^p M/U$$
	Thus $\prod^p M/U$ is not closed under addition.
	\item \label{badex-N} Let $L$ be the two-sorted language $\seq{N_1, N_2; +_1, \times_1, 1_1; +_2, \times_2, 1_2; \times_{1, 2}}$ where $\times_{1, 2}:N_1\times N_2 \to N_1$.  Take $M = \seq{\mathbb{N}, \mathbb{N}'; +, \times, 1; +', \times', 1'; \times^*}$ where $\mathbb{N}$ and $\mathbb{N}'$ are disjoint copies of the naturals and $\times^*$ is also normal multiplication.  Then this structure omits the type of a nonstandard element of the second sort $p(x) = \{ N_2(x) \wedge (1 + \dots + 1 \neq x) : n < \omega\}$.  Then $\prod^p M / U$ is a structure.  In particular, $N_2$ remains standard but $N_1$ is just $\prod \mathbb{N}/U$.  To see the failures of \L o\'{s}' Theorem described above, 
\begin{itemize}
	\item the formula $\psi(x)\equiv ``\exists y \in N_2 (1_1 \times_{1, 2} y = x)$ is true of all $n \in N_1^M$, but is not true of $[n \mapsto n]_U \in N_1^{\prod^p M/U}$.
	\item the sentence $\phi \equiv \te{``}\forall x \in N_1 \exists y \in N_2 ( 1_1 \times_{1, 2} y = x)\te{''}$ is true in $M$, but not in $\prod^p M/U$ for the above reason.
\end{itemize}
\end{enumerate}
\end{example}

Note that the first example shows that the class of classically valued fields does not fit into the framework described here;  compactness results in that class will be explored in Boney \cite{cacvf}.

These examples and Example \ref{badpower-ex} below give an indication of when things don't fit nicely into this framework.  Section \ref{example-sec} collects several positive examples.

\subsection{Structure}

For $\prod^\Gamma M_i/U$ to be a structure, all that is necessary is that $\prod^\Gamma M_i/U$ is closed under functions.  This means that, if $[f_0,]_U \dots, [f_{n-1}]_U \in \prod^\Gamma M_i/U$ have witnesses to their inclusion and $F$ is a function of $\prod M_i/U$, then $F\left([f_0]_U, \dots, [f_{n-1}]_U\right)$ has a witness as well.  However, in many cases, there is a degree of uniformity where tuples with the same sequence of witnesses are always mapped to an element with a fixed witness.

\begin{defin} \label{gamclose-def}
$\{M_i : i \in I\}$ is $\Gamma$-closed iff for all $n$-ary functions $F$ of $L$, there is a function $g_F$ that takes in $n$ choice functions on $\Gamma$ and outputs a choice function on $\Gamma$ such that
\begin{center}
for all $[f_1]_U, \dots, [f_{n}]_U \in \prod^\Gamma M_i/U$ with witnesses $\mathcal{C}_1, \dots, \mathcal{C}_{n}$ , we have that, for each $i \in I$ and $p \in \Gamma$,
$$M_i \vDash \neg \phi\left(F\left(f_1(i), \dots, f_{n}(i)\right)\right)$$
where $\phi = g_F\left(\mathcal{C}_1, \dots, \mathcal{C}_{n}\right)(p)$
\end{center}
\end{defin}

Abelian torsion groups are an example of this:  the order of $h + k$ can be computed from the orders of $h$ and $k$; see Section \ref{atg-ex}.  There, $g_+$ takes in natural numbers (representing choice functions on the singleton set of the torsion type) and outputs a natural number such that $g_+(o(h), o(k))$ is an order of $h+k$. \\
It is clear that if $\{M_i : i \in I\}$ is $\Gamma$-closed, then $\prod^\Gamma M_i/U$ is a structure for all ultrafilters $U$; the witness to $F\left([f_0]_U, \dots, [f_{n-1}]_U\right)$ is $g_F\left(\mathcal{C}_1, \dots, \mathcal{C}_{n}\right)$.\\

The main advantage of $\Gamma$-closedness is in the study of classes of models omitting $\Gamma$ when $\Gamma$ is finite because this property is captured by the first order theory of $M$.  We say that a class $EC(T, \Gamma)$ is $\Gamma$-closed iff every collection of models from it is $\Gamma$-closed.

\begin{prop} \label{gamclofo-prop}
Suppose $\Gamma$ is finite and $T$ is $\forall(\Gamma \cup \neg\Gamma)$-complete.  Then $EC(T, \Gamma)$ is $\Gamma$-closed iff some collection $\{ M_i \in EC(T, \Gamma) : i \in I\}$ is.
\end{prop}

The notation ``$\forall(\Gamma \cup \neg\Gamma)$-complete'' means that $T$ decides all first order sentences whose quantifiers are a universal followed by a quantifier string that appears in $\Gamma$ or a universal followed by the negation of such a string; similar expressions have the obvious meaning.

{\bf Proof:}  Being $\Gamma$-closed can be expressed by the following scheme: for each $F \in L$, sequence of choice functions $\mathcal{C}_0, \dots, \mathcal{C}_{n-1}$, and $p \in \Gamma$, include the sentence
$$\forall x_0, \dots, x_{n-1} \left( \left(\bigwedge_{q \in \Gamma, i<n} \neg \mathcal{C}_i\left(q\right)(x_i) \right) \to \neg g_F\left( \mathcal{C}_0, \dots, \mathcal{C}_{n-1} \right) (p) \left( F(x_0, \dots, x_{n-1}) \right)\right)$$
where $g_F$ is the witness for $F$ from the definition of $\Gamma$-closed.  Since $\Gamma$ is finite, this is first order of the desired complexity. \hfill \dag\\

If $\prod^\Gamma M_i/U$ is a structure, then it is a substructure of $\prod M_i/U$.  This immediately gives us a universal version of \L o\'{s}' Theorem.

\begin{theorem}[Universal \L o\'{s}' Theorem]\label{weaklos}
Suppose $\prod^\Gamma M_i/U$ is a structure. If $\phi(x_0, \dots, x_n)$ is a universal formula and $[f_0]_U, \dots, [f_{n-1}]_U \in \prod^\Gamma M_i/U$, then
$$\{i \in I : M_i \models \phi(f_0(i), \dots, f_{n-1}(i)) \} \in U \implies \prod^\Gamma M_i/U \models \phi([f_0]_U, \dots, [f_{n-1}]_U)$$
\end{theorem}

{\bf Proof:} The key point is that, since $\prod^\Gamma M_i/U$ is a structure, it is a substructure of the full ultraproduct $\prod M_i/U$.  Thus, universal formulas transfer from $\prod M_i/U$ to $\prod^\Gamma M_i/U$.  \comment{Let $\phi(\bx)$ be universal and $[f_0]_U, \dots, [f_{n-1}]_U \in \prod^\Gamma M_i/U$.  Then
\begin{eqnarray*}
\{i \in I: M_i \vDash \phi(f_0(i), \dots, f_{n-1}(i) \} \in U &\implies& \prod M_i /U \vDash \phi([f_0]_U, \dots, [f_{n-1}]_U) \\
&\implies& \prod^\Gamma M_i /U \vDash \phi([f_0]_U, \dots, [f_{n-1}]_U)
\end{eqnarray*}}
\hfill \dag \\

\begin{remark} \label{inductlos-remark}
A proof by induction on formula complexity (mirroring the proof of the standard version of \L o\'{s}' Theorem) is also possible.  This proof is longer, but provides extra information: if $\phi$ and $\psi$ are formulas that transfer from a $U$-large set of $M_i$ to $\prod^\Gamma M_i/U$, then conjunction, disjunction, and universal quantification preserves this transfer, while negation reverses it.  This finer analysis is used in Section \ref{torcomp-sec}.
\end{remark} 

This has important implications for $\Gamma$ consisting of existential types.

\begin{prop} \label{comthe}
Suppose $\prod^\Gamma M_i/U$ is a structure and the types of $\Gamma$ contain only existential formulas.
\begin{enumerate}
	\item $\prod^\Gamma M_i/U$ omits $\Gamma$.
	
	\item If $\Gamma$ is finite and each $M_i$ satisfies a common $\exists\forall$ theory $T_{\exists\forall}$, then $\prod^\Gamma M_i/U \vDash T_{\exists \forall}$.

\end{enumerate}
\end{prop}

{\bf Proof:}\begin{enumerate}
	\item  Let $[f]_U \in \prod^\Gamma M_i/U$ and $p \in \Gamma$.  By definition, there is $\phi_p\in p$ such that $\{i \in I: M_i \vDash \neg \phi_p\left(f(i)\right)\}\in U$.  Since $\neg \phi_p$ is universal, $\prod^\Gamma M_i/U \vDash \neg \phi_p(m)$ by Universal \L o\'{s}' Theorem \ref{weaklos}. 
	\item Let $\exists \bx \psi(\bx)$ be in $T_{\exists\forall}$ with $\psi$ universal and let $i_0 \in I$.  Then there are $m_1, \dots, m_n \in M_{i_0}$ such that $M_{i_0} \vDash \psi(m_1, \dots, m_n)$.  Because $M_{i_0}$ omits $\Gamma$, for each $p \in \Gamma$ and $\ell=1, \dots, n$, there is $\phi^\ell_p \in p$ such that $M_{i_0} \vDash \neg \phi^\ell_p(m_\ell)$.  Then, 
	$$M_{i_0} \vDash \exists \bx \left( \psi(x_1, \dots, x_n) \wedge \bigwedge_{p \in \Gamma; \ell \leq n} \neg \phi_p^\ell(x_\ell) \right)$$
	This is part of $T_{\exists\forall}$, so for each $i \in I$, there is $m_1^i, \dots, m_1^n \in M_i$ such that
	$$M_i \vDash \psi(m_1^i, \dots, m_n^i) \wedge \bigwedge_{p\in\Gamma; \ell \leq n} \neg \phi^\ell_p(m_\ell^i)$$
	Thus $[g_\ell]_U \in \prod^\Gamma M_i / U$, where $g_\ell(i) = m_\ell^i$, and, by Universal \L o\'{s}' Theorem \ref{weaklos}, we have that
	$$\prod^\Gamma M_i / U \vDash \psi\left( [g_1]_U, \dots, [g_n]_U\right)$$
\end{enumerate}\hfill \dag\\

Indeed, if $\Gamma$ is finite (but not necessarily existential), a similar proof shows that if each $M_i$ satisfy a common $\exists (\neg \Gamma \cup \forall)$-theory, then $\prod^\Gamma M_i/U$ models the $\exists\forall$ part of the common theory.  Note that the $\exists\forall$ level is sharp as Example \ref{badex}.(2) gives an example of a $\Gamma$-ultrapower that doesn't have the same $\forall\exists$ theory.

\subsection{Elementary Substructure} \label{elem-subsec}

For this subsection, assume that $\prod^\Gamma M_i/U$ is a structure.

    The second line of analysis of the $\Gamma$-ultraproduct is finding the class of formulas $\phi(\bx)$ such that for all $[f_0]_U, \dots, [f_{n-1}]_U \in \prod^\Gamma M_i/U$,
$$ \{ i \in I : M_i \vDash \phi\left(f_0(i), \dots, f_{n-1}(i) \right) \} \in U \iff \prod^\Gamma M_i/U \vDash \phi\left([f_0]_U, \dots, [f_{n-1}]_U\right)$$
We know that this class contains the universal formulas and, indeed, is closed under universal quantification.  However, Example \ref{badex}.(2) shows that existential quantification causes problems.  The problem is that the witnesses to an existential formula involving parameters that omit $\Gamma$ uniformly might not omit $\Gamma$ uniformly.   

When $\Gamma$ is finite, some level of existential quantification is allowed by essentially forcing a witness to exist as part of the condition on the existential.

\begin{observation}\label{doublestar}
Suppose $\Gamma$ is finite.  If \L o\'{s}' Theorem holds for $\phi(x, \by)$ then it also holds for
$$\exists x \left( \phi(x, \by) \wedge \bigwedge_{p \in \Gamma} \neg \mathcal{C}(p)(x) \right)$$
for any choice function $\mathcal{C}$ on $\Gamma$.  Recalling Ben Yaacov's work on metric ultraproducts in unbounded metric structures, this condition is similar to his requirement that the formula is bounded \cite[Definition 2.7]{byunbounded}.
\end{observation}

In general, there are two main ways of guaranteeing the transfer of all existential statements: $\Gamma$-niceness and quantifier elimination.\\

$\Gamma$-niceness is the appropriate generalization of $\Gamma$-closed to the situation of existentials.

\begin{defin} \label{gamnice-def}
$\{M_i : i \in \Gamma\}$ is $\Gamma$-nice iff for all existential formulas $\psi :\equiv \exists x \phi(x, \by)$ of $L$, there is a function $g_\psi$ that takes in $\ell(\by)$ choice functions on $\Gamma$ and outputs a choice function on $\Gamma$ such that 
\begin{center}
for all $[f_1]_U, \dots, [f_n]_U \in \prod^\Gamma M_i / U$ with witnesses $\mathcal{C}_1, \dots, \mathcal{C}_n$ and $i \in I$, if $M_i \vDash \exists x \phi\left(x, f_1(i), \dots, f_n(i)\right)$, then there is $m \in M_i$ such that, for all $p \in \Gamma$,
$$M_i \vDash \phi\left(m, f_1(i), \dots, f_n(i)\right) \wedge \neg \chi(m)$$
where $\chi = g_\psi\left(\mathcal{C}_1, \dots, \mathcal{C}_n\right) (p)$.
\end{center}
\end{defin}

The following basic facts about $\Gamma$-niceness are obvious.  Recall from the introduction that the data refers to the collection of the language, the types to be omitted, and the structures that we wish to take the $\Gamma$-ultraproduct of.

\begin{prop} \label{nicefacts}
\begin{enumerate}

	\item If the data is $\Gamma$-nice, then it is $\Gamma$-closed.
	
	\item The data is $\Gamma$-nice iff there is a skolemization of the data that is $\Gamma$-closed.
	
	\item If $\Gamma$ is finite, then being $\Gamma$-nice is first-order expressible.

\end{enumerate}
\end{prop}

{\bf Proof:} For (1), set $g_F:=g_{\exists x (F(\by) = x)}$.  For (2), take $g_{F_{\exists x \phi(x; \by)}} = g_{\exists x \phi(x; \by)}$.  For (3), the proof follows as in Proposition \ref{comthe}.(2).\hfill\dag\\

The main use of $\Gamma$-niceness is as a sufficient condition for \L o\'{s}' Theorem to hold.

\begin{theorem}[\L o\'{s}' Theorem] \label{lostheorem}
Suppose the data is $\Gamma$-nice and $U$ is an ultrafilter on $I$. If $\phi(x_1, \dots, x_n)$ is a formula and $[f_1]_U, \dots, [f_{n}]_U \in \prod^\Gamma M_i/U$, then
$$\{i \in I : M_i \models \phi\left(f_1(i), \dots, f_{n}(i)\right) \} \in U \iff \prod^\Gamma M_i/U \models \phi\left([f_1]_U, \dots, [f_{n}]_U\right)$$
\end{theorem}

{\bf Proof:} By Proposition \ref{nicefacts} and Theorem \ref{weaklos}, all that needs to be shown is that adding an existential quantifier maintains transfer from ``true in $U$-many $M_i$'s'' to ``true in $\prod^\Gamma M_i/U$.''  That is, suppose $\phi(\by) = \exists x \psi(x, \by)$ such that, for all $[f_0]_U, \dots, [f_n]_U \in \prod^\Gamma M_i/U$, 
	$$\{i \in I : M_i \models \psi\left(f_0(i), \dots, f_{n}(i)\right) \} \in U \iff \prod^\Gamma M_i/U \models \psi\left([f_0]_U, \dots, [f_{n}]_U\right)$$

We want to show that, for all $[f_1]_U, \dots, [f_n]_U \in \prod^\Gamma M_i/U$, 
	$$ X:=\{i \in I : M_i \models \exists x\psi(x,f_1(i), \dots, f_{n-1}(i)) \} \in U \implies \prod^\Gamma M_i/U \models \exists x  \psi(x, [f_1]_U, \dots, [f_{n-1}]_U)$$
Suppose we have such a tuple with witnesses $\mathcal{C}_1, \dots, \mathcal{C}_n$.  By $\Gamma$-niceness, for each $i \in X$, there is $m_i \in M_i$ as in the definition.  Then $[i \mapsto m_i]_U$ is in $\prod^\Gamma M_i/U$, as witnessed by $g_\psi(\mathcal{C}_1, \dots, \mathcal{C}_n)$ and $\{i \in I : M_i \models \psi\left(m_i, f_1(i), \dots, f_n(i)\right) \} \in U$.  By the induction assumption,
$$\prod^\Gamma M_i / U \vDash \psi\left([i\mapsto m_i]_U, [f_1]_U, \dots, [f_n]_U \right)$$
as desired. \hfill \dag\\

Once we have the full strength of \L o\'{s}' Theorem, we are guaranteed the resulting structure omits the desired types. 

\begin{prop}[Type Omission]\label{typeomission}
Suppose the data is $\Gamma$-nice (or just \L os' Theorem holds).  Then $\prod^\Gamma M_i/U$ omits each type in $\Gamma$.
\end{prop}

{\bf Proof:} Let $[f]_U \in \prod^\Gamma M_i/U$.  This is witnessed by some $\mathcal{C}$.  For each $p \in \Gamma$,
$$\{i \in I:M_i \vDash \neg \mathcal{C}(p)\left(f(i)\right) \} \in U$$
By Theorem \ref{lostheorem}, 
$$\prod^\Gamma M_i/U \vDash \neg \mathcal{C}(p)\left([f]_U\right)$$
So every element of $\prod^\Gamma M_i /U$ does not realize any type from $\Gamma$. \hfill \dag\\

Summarizing our results so far, we have the following.

\begin{cor}
If the data is $\Gamma$-nice, then $\prod^\Gamma M_i/U$ is an $L$-structure that satisfies \L o\'{s}' Theorem and omits every type in $\Gamma$.  In particular, if $M_i \in EC(T, \Gamma)$ for all $i \in I$, then $\prod^\Gamma M_i/U \in EC(T, \Gamma)$.
\end{cor}

{\bf Proof:} By Theorems \ref{lostheorem} and \ref{typeomission}.\hfill \dag\\

We now turn to another method for proving \L o\'{s}' Theorem that seems more ad-hoc, but has proven more useful in practice: quantifier elimination.  This method involves directly proving that the $\Gamma$-ultraproduct is a structure that models the theory $T$ and, if $T$ has only partial quantifier elimination, proving that \L o\'{s}' Theorem holds for the necessary class of formulas.  Examples of this are $DLOG\mathbb{Z}$ (Section \ref{dlogz-ex})
\comment{, algebraic $DCF$\footnotei{WB:Name} (Section \ref{dcf-ex}),}
 and torsion modules over PIDs (Section \ref{tormod-sec}).  This final example makes use of the full generality of the following proposition since modules only have quantifier elimination to p. p. formulas.

\begin{prop} \label{qelos-prop}
Suppose each $M_i \vDash T$, $\prod^\Gamma M_i/U \vDash T$, $T$ has quantifier elimination to $\Delta$-formulas, and \L o\'{s}' Theorem for $\Delta$-formulas holds.  Then the full \L o\'{s}' Theorem holds.
\end{prop}
 
{\bf Proof:} Immediate. \hfill \dag\\

\subsection{Ultrapowers} \label{ultra-subsec}
We now turn our attention to ultrapowers, where $M_i = M$ for all $i \in I$.  In this case, set $\j:M \to \prod^\Gamma M_i/U$ to be the ultrapower map by $\j(m) = [i \mapsto m]_U$.  This function is well-defined even if $\prod^\Gamma M_i/U$ is not a structure, and the statement of \L os' Theorem is equivalent to $\j$ being an elementary embedding.  We would also like to know when this construction gives rise to a proper extension.  Unfortunately, this is not always the case.

\begin{example} \label{badpower-ex}
Let $U$ be an ultrafilter on $I$.  Take the ultraproduct of $\mathbb{N} = \seq{\omega, +, \cdot, <}$ omitting $p(x) = \{ x > n : n < \omega \}$; then $\j:\mathbb{N} \cong \prod^p \mathbb{N}/U$.
\end{example}

{\bf Proof:}  The data is obviously $\Gamma$-closed: $g_+(x, y) = x + y$ and $g_\cdot(x, y) = x \cdot y$.  This is enough to make the conclusion well-formed, i. e. $\prod^p \mathbb{N}/U$ is a structure.  If $f \in \prod^p \mathbb{N}$, then there is some $k_f < \omega$ such that $f(i) < k_f$ for all $i \in I$.  Since $U$ is $\omega$-complete (as are all ultrafilters) and $k_f$ is finite, there is some $n_f < k_f$ such that $\{ i \in I : f(i) = n_f\} \in U$.  Thus, $[f]_U = [ i \mapsto n_f]_U$ and the mapping $h: \prod^p \mathbb{N}/U \to \mathbb{N}$ by $h([f]_U) = n_f$ is an isomorphism. \hfill \dag\\

This did not give rise to a new model because the choice function witnessing $f \in \prod^p \mathbb{N}$, here characterized by a single natural number, determined which element of $\mathbb{N}$ the function represented.  In order to ensure that $\j$ is not surjective, we need to ensure that there are many choices that give rise to the same $\mathcal{C}$.  Indeed, this characterization is reversible.

\begin{theorem} \label{propext}
Suppose $M$ omits $\Gamma$.  Then the following are equivalent:
\begin{enumerate}

	\item There is infinite $X \subset M$ and choice function $\mathcal{C}$ such that, for all $m \in X$ and $p \in \Gamma$, 
	$$M \vDash \neg \mathcal{C}(p) \left(m\right)$$
	
	\item There is a nonprincipal ultrafilter $U$ such that $\j : M \to \prod^\Gamma M/U$ is not surjective.
	
\end{enumerate}
\end{theorem}	
	
{\bf Proof:} For (1) implies (2), let $U$ be any nonprincipal ultrafilter on $\omega$.  Then let $f:\omega \to X$ enumerate distinct members of $X$.  By definition of $X$, $f \in \prod^\Gamma M$.  Since $U$ is nonprincipal, $f$ is not $U$-equal to any constant function.  Thus, $[f]_U$ is an extra element in $\prod^\Gamma M/U$.	\\
For (2) implies (1), let $[f]_U$ be a new element of $\prod^\Gamma M/U$.  Then, for each $m \in M$, 
$$X_f \cap \{i\in I : f(i) = m \} \not\in U$$
Then $f''X_f \subset M$ is infinite and, taking $\mathcal{C}$ to be the choice function witnessing $f \in \prod^\Gamma M$, we have that, for each $m \in f''X_f$ and $p \in \Gamma$,
$$M \vDash \neg \mathcal{C}(p) \left(m\right)$$
\hfill \dag

\begin{cor} \label{upext-cor}
If $\|M\| > \prod_{p \in \Gamma} |p|$, then there is an ultrafilter $U$ such that $\j:M \to \prod^\Gamma M/U$ is not surjective.
\end{cor}

{\bf Proof:} For each $x \in M$, pick a choice function $\mathcal{C}^x$ such that, for all $p \in \Gamma$, $M \vDash \neg \mathcal{C}^x(p)(x)$.  There are $\prod |p|$ many possible values for $\mathcal{C}^x$.  Since $\|M\|$ is greater than this, there must be some infinite $X \subset |M|$ such that the choice is constant.  Then apply Theorem \ref{propext}.\hfill \dag\\

\begin{cor}
Suppose $p$ is countable and $M$ is uncountable.  If $U$ is nonprincipal, then $\prod^p M/U \not\cong M$.
\end{cor}

\subsection{Changing the language}

One strength of the ultraproduct is its robustness under changing the language.  Unfortunately, the $\Gamma$-ultraproduct does not share this robustness.  However, some results remain, which gives rise to a form of $\Gamma$-compactness in Theorem \ref{gamcomp}, and there are sometimes natural conditions, such as in Section \ref{bs-ex}, that determine when the expansions are well-behaved. 

As with quantifier elimination, adding constants does not impact the properties discussed above ($\Gamma$-closed, etc.).  As an example, we show that $\Gamma$-niceness is preserved.

\begin{prop}\label{constexp-prop}
Suppose $\{M_i : i \in I\}$ is $\Gamma$-nice and $M_i^*$ is an expansion of $M_i$ by constants $\{c_j : j<\kappa\}$ such that, for all $j < \kappa$, $c_j^{M_i}$ omits each type at the same place for all $i \in I$.  Then $\{M_i^* : i \in I\}$ is $\Gamma$-nice.
\end{prop}

{\bf Proof:} Let $\psi^*(\by) = \exists x \phi^*(x, \by)$ be an existential formula in the expanded language.  Then, there are new constants $\bc$ such that $\exists x \phi(x, \by, \bc) = \psi^*$ and $\phi$ is a formula in the original language.  Set $g_{\psi^*}(\bz) : = g_{\psi}(\bz, \bar{\mathcal{C}})$, where $\bar{\mathcal{C}}$ are the choice functions for the $\bc$; these exist by hypothesis.  Then $g_{\psi^*}$ witness the $\Gamma$-niceness.\hfill \dag\\

In general, expanding the language by functions or relations does not preserve these properties.  For instance, an added function might pick out elements that omit the types ``wildly'' on a domain that omits the types at the same place.  This is of course unfortunate because these are the expansions that are most often useful.  However, this means that the study of when expanding the language preserves these properties is of great interest.  One example is given in Subsection \ref{bs-ex}.  Another example, given below, shows that this investigation allows us to get compactness results outside of omitting types classes.

Set $T$ to be the $L(Q)$-theory that says $E$ is an equivalence relations and each class is countable;  this is a first order axiom and $\forall x \neg Q y E(x, y)$.  This is the most basic example of a quasiminimal class and of a non-finitary AEC (the strong substructure is substructure plus equivalence classes don't grow; see Kirby \cite{qmc} for an explicit description of this class and overview of quasiminimal classes).  Thus, this is not a type-omitting class, but there is a well-known method that allows the expression of $L_{\omega_1, \omega}(Q)$ in terms of $L_{\omega_1, \omega}$ in an expanded language (see, for instance, the proof of \cite[Theorem 5.1.8]{baldwinbook}).  This expansion is not canonical and typically gives rise to non-unary types.  However, in this example, the combination of the facts that the $L(Q)$-subformula has only one free variable and the fact that the quasiminimal closure is trivial allows us to get a result.

Expand $L$ by adding countably many unary predicates $\{R_n(x):n<\omega\}$ and expand a model of $T$ by making $R_n$ true of exactly one member of each equivalence class.  Set $T^*$ to be the first order part of $T$ plus $\forall x \exists ! x\left(R_n(y) \wedge E(x, y)\right)$ for each $n < \omega$ and set $p = \{\neg R_n(x) : n < \omega\}$.  The following is straightforward.

\begin{claim}
Let $\{M_i : i \in I\}$ be models of $T$ and $U$ an ultrafilter on $I$.  If $M_i^*$ is an expansion of $M_i$ to a model of $T^*$ omitting $p$, then 
$$\left(\prod^p M_i^*/U\right) \rest L \vDash T$$
and \L o\'{s}' Theorem holds.  Moreover, $\left(\prod^p M_i^*/U\right) \rest L$ does not depend on the choice of the expansion.
\end{claim}

This is a very basic example and the consequences are more easily obtained by analyzing it as a quasiminimal class.  However, it gives hope that more intractable $L(Q)$ classes can be analyzed via the $\Gamma$-ultraproduct.

\subsection{A different approach} \label{diffapp-sec}

We have given a description of this construction that is intimately tied to the ultraproduct.  However, if $\Gamma$ is a finite set of types, then there is an equivalent way of constructing the $\Gamma$-ultraproduct.

Given a model $M$, set the $\Gamma$-hull of $M$ to be 
$$\Gamma(M):=\{m \in M: \forall p \in \Gamma, m \text{ does not realize }p\}$$
If $\Gamma$ is finite, then $\Gamma\left(\prod M_i/U\right) = \prod^\Gamma M_i/U$.  In the general study of $\Gamma$-hulls, i. e., not just when taking the $\Gamma$-hull of an ultraproduct, many of the same issues arise in the analysis of the $\Gamma$-hull as the $\Gamma$-ultraproduct, but with the relationship between $\Gamma(M)$ and $M$ taking center stage.  For instance, Section \ref{elem-subsec} would be replaced by an exploration of how elementary $\Gamma(M)$ is in $M$.

This paper focused on the $\Gamma$-ultraproduct over the $\Gamma$-hull for two related reasons.  First, the main goal of viewing nonelementary classes as averageable is that some form of compactness holds there.  Thus, working with ultraproducts is very natural.  Second, some classes have a stronger \L o\'{s}' Theorem between $\{M_i : i \in I\}$ and $\prod^\Gamma M_i/U$ than there is (in general) elementarity between $\Gamma(M)$ and $M$.  This again makes $\prod^\Gamma M_i/U$ the natural choice.  An example of the second is abelian torsion groups.  For instance, $tor\left(\mathbb{Z} \oplus \mathbb{Z}_2\right) = \mathbb{Z}_2$ and $\mathbb{Z} \oplus \mathbb{Z}_2$ are very different, but the full \L o\'{s}' Theorem holds for this class by Theorem \ref{fulllos-thm}.

In general, we have, as sets,
$$\prod^\Gamma M_i / U \subset \Gamma\left(\prod M_i/U\right) \subset \prod M_i/U$$
It would be interseting to find a nonelementary class (with $\Gamma$ necessarily infinite) where the $\Gamma$-hull of the ultraproduct was the proper structure to analyze, e. g., it is different from the $\Gamma$-ultraproduct and the $\Gamma$-hull is in the class, but the $\Gamma$-ultraproduct is not.

\section{Averageable Classes}\label{avclasses} \label{aecsection}

We now consider classes that are well behaved under some $\Gamma$-ultraproduct.  We use the language of Abstract Elementary Classes (AECs) here; Baldwin \cite{baldwinbook}, Grossberg \cite{g-bilgi}, and Shelah \cite{shelahaecbook} are the standard references.  The main object of study are classes of models that omit all types from $\Gamma$ such that \emph{all} sets of models from there are $\Gamma$-closed or $\Gamma$-nice and that the strong substructure relation under consideration is preserved under ultraproducts.  The examples we consider all fall into one of two cases:
\begin{enumerate}
	\item $(EC(T, \Gamma), \subset)$ when $T$ is $\exists\forall$, the types of $\Gamma$ are existential, and $EC(T, \Gamma)$ is $\Gamma$-closed.
	\item $(EC(T, \Gamma), \prec)$ when $EC(T, \Gamma)$ is $\Gamma$-nice.
\end{enumerate}
The reader can easily focus on these, but we introduce a joint generalization of these cases to keep from stating our results twice and for potential future applications.

\begin{defin}
An AEC $K$ is averageable iff there is a collection of first-order formulas $\eff$ that contains all atomic formulas and a collection of unary types $\Gamma$ such that
\begin{itemize}
	\item for all $\phi(x) \in p \in \Gamma$, $\neg \phi(x) \in \eff$;
	\item the strong substructure $\prec_K$ is $\eff$-elementary substructure;
	\item given $\{M_i \in K : i \in I\}$ and an ultrafilter $U$ on $I$, $\prod^\Gamma M_i / U \in K$ and this $\Gamma$-ultraproduct satisfies \L o\'{s}' Theorem for the formulas in $\eff$; and
	\item each $M \in K$ omits each $p \in \Gamma$.
\end{itemize}
\end{defin}

Being an averageable class gives a very strong compactness result within the incompact framework of AECs; this was seen using large cardinals in Boney \cite{tamelc} and exploited in Boney and Grossberg \cite{bgcoheir} axiomatically.

The first property of averageable classes is a (much) better Hanf number.

\begin{prop}\label{nomax}
Suppose that $K$ is averageable and set $\kappa = \prod |p|$.  Then $K_{> \kappa}$ has no maximal models.
\end{prop}

This follows directly from Corollary \ref{upext-cor}.  The normal Hanf number is $\beth_{(2^{|T|})^+}$, while $\prod |p| \leq 2^{|T|}$.

A stronger result builds on Proposition \ref{constexp-prop}.  Because they omit types, averageable classes are not compact.  However, they satisfy a strong approximation to compactness that we call local compactness.

\begin{theorem}[Local Compactness]\label{gamcomp}
Suppose that $\left(EC(T, \Gamma), \prec\right)$ is an averageable class with $\mathcal{F}$; $T^*$ is an extension of $T$ in $\tau^*:= \tau(T) \cup \{c_i : i < \kappa\}$ such that each sentence comes from substituting new constants into a formula from $\mathcal{F}$; and choice functions $\{ \mathcal{C}_i : i < \kappa\}$ on $\Gamma$ with $T_0 = \{ \neg \mathcal{C}_i(p)\left(c_i\right) : i < \kappa, p \in \Gamma\}$.

If, for every finite $T^- \subset T^*$, there is $M^* \vDash T^- \cup T_0$ that omits all of $\Gamma$, then $T^*$ has a model omitting all of $\Gamma$.
\end{theorem}

Looking at the space $gS^n(\emptyset)$ of syntactic $n$-types that are realized in some model of $EC(T, \Gamma)$ when $\Gamma$ is finite leads to the name local compactness.  Equipped with the standard topology, this space is not compact, but it is locally compact.  If $p \in gS^n(\emptyset)$ and $\mathcal{C}_1, \dots, \mathcal{C}_n$ are the corresponding choice functions, then 

$$\left\llbracket \bigwedge_{p \in \Gamma, i < n} \neg \mathcal{C}_i(p)(\left(x_i\right) \right\rrbracket$$

is a compact clopen neighborhood of $p$ (the compactness follows by our theorem).  When $n$ or $\Gamma$ is infinite, this set is no longer clopen (so the space is no longer locally compact), but the set above is still compact and typically contains more than just $p$.

{\bf Proof:}  For each finite $T^- \subset T^*$, let $M^*_{T^-}$ be the model advertised in the hypothesis and let $M_{T^-}$ be its restriction to $\tau(T)$.  Let $U$ be a fine ultrafilter on $P_\omega T^*$ and set 
$$M := \prod^\Gamma_{T^- \in P_\omega T^*} M_{T^-}/U$$
By averageability, $M \in EC(T, \Gamma)$ and \L o\'{s}' Theorem holds for $\mathcal{F}$.  Now expand $M$ to a $\tau^*$ structure $M^*$ by setting $c_i^{M^*} := [T^- \mapsto c_i^{M^*_{T^-}}]_U$.  Note that $\mathcal{C}_i$ is the witness to $c_i^{M^*} \in M$, so this is a valid definition.  Now we claim that $M^*$ is the model satisfying $T^*$ and omitting $\Gamma$.  Since $M \in EC(T, \Gamma)$, it satisfies $T$ and omits $\Gamma$ and naming additional constants does not change this.  Now let $\phi \in T^*$ be a new sentence.  Then it is of the form $\psi(c_{i_1}, \dots, c_{i_n})$ for some $\psi \in \mathcal{F}$.  By the fineness of the ultrafilter,
$$[\phi] = \{ T^- \in P_\omega T^* : M_{T^-} \vDash \psi( c^{M^*_{T^-}}_{i_1}, \dots, c^{M^*_{T^-}}_{i_n}) \} \in U$$
Since \L o\'{s}' Theorem holds for $\psi$, this means that 
$$M^* \vDash \psi\left( c_{i_1}^{M^*}, \dots, c_{i_n}^{M^*}\right)$$
Since this holds for every $\phi \in T^*$, we have $M^* \vDash T^*$, as desired. \hfill \dag\\

Although complex to parse, local compactness has a very nice corollary.

\begin{cor}
Suppose $EC(T, \Gamma)$ is $\Gamma$-nice and $p$ is a type such that every finite subset is a realizable in a model of $EC(T, \Gamma)$ with some fixed witness.  Then $p$ is realized in a model of $EC(T, \Gamma)$.
\end{cor}

Many other uses of compactness follow similarly, such as a criteria for amalgamation similar to first-order (see \cite[Theorem 6.5.1]{hodges}).

\begin{cor}
Suppose $EC(T, \Gamma)$ is $\Gamma$-nice.  Then it has amalgamation.
\end{cor}

Indeed, much of the rest of this section is more easily proven with the local compactness result above.  However, we have lost some generality by restricting ourselves to averageable classes of the form $EC(T, \Gamma)$.  Although this is the form of our examples, we prove the following results in greater generality and, therefore, without local compactness.

Following \cite{tamelc}, we get the following result about tameness and type shortness.  Tameness and type shortness are locality results for Galois types.  We only prove type shortness (and not the tameness follows), so we define that here.  For full definitions, see, e. g., \cite[Section 3]{tamelc}.

\begin{defin}
Let $K$ be an AEC and $I$ a linear order.
\begin{itemize}
	\item Given $M \prec_K N_1, N_2$ from $K$ and $\seq{x_i^\ell \mid i \in I}$, we say that
	$$\left( \seq{x_i^1 \mid i \in I}, M, N_1\right) E_{AT} \left( \seq{x_i^2 \mid i \in I}, M, N_2\right)$$
	iff there is $N^* \in K$ and $g_\ell : N_\ell \to_M N^*$ such that $g_1(x_i^1) = g_2(x_i^2)$ for all $i \in I$.
	\item Given $M \prec_K N_1, N_2$ from $K$ and $\seq{x_i^\ell \mid i \in I}$, we say that $\left( \seq{x_i^1 \mid i \in I}, M, N_1\right)$ and $\left( \seq{x_i^2 \mid i \in I}, M, N_2\right)$ have the same Galois type, written
	$$gtp(\seq{x_i^1 \mid i \in I}/ M; N_1)=gtp(\seq{x_i^2 \mid i \in I}/ M; N_2)$$
	iff they are related by the transitive closure of $E_{AT}$.  The length of the Galois type $gtp(\seq{x_i^1 \mid i \in I}/ M; N_1)$ is the index $I$
	\item $K$ is fully $<\kappa$-type short iff for all $I$ and for all Galois types $gtp(\seq{x_i^1 \mid i \in I}/ M; N_1)$ and $gtp(\seq{x_i^2 \mid i \in I}/ M; N_2)$, we have 
	$$gtp(\seq{x_i^1 \mid i \in I}/ M; N_1)=gtp(\seq{x_i^2 \mid i \in I}/ M; N_2)$$
	iff for all $I_0 \subset I$ of size $< \kappa$
	$$gtp(\seq{x_i^1 \mid i \in I_0}/ M; N_1)=gtp(\seq{x_i^2 \mid i \in I_0}/ M; N_2)$$
\end{itemize}
\end{defin}

Note that, if $K$ satisfies amalgamation, then $E_{AT}$ is already transitive.  Full $<\omega$-type shortness follows from the assertion that Galois types are syntactic (in some sublogic of $L_{\infty, \omega}$).  For examples of AECs that are \emph{not} type short, see Baldwin-Shelah \cite{nonlocality}.

\begin{theorem}\label{tameresult}
Suppose $K$ is averageable.  Then $K$ is fully $<\omega$-tame and -type short.
\end{theorem}

This relies on the following lemma, which says that the ultraproduct of $K$-embeddings is also a $K$-embedding.

\begin{lemma}\label{blank}
Suppose that $\seq{M_i : i \in I}$ and $\seq{N_i : i \in I}$ and $f_i:M_i \to N_i$ is a $K$-embedding.  Then $f:\prod^\Gamma M_i / U \to \prod^\Gamma N_i / U$ by $f([i \mapsto m_i]_U) = [i \mapsto f_i(m_i)]_U$ is a $K$-embedding.
\end{lemma}

{\bf Proof:} First, we need to know that $[i\mapsto f_i(m_i)]_U$ is in $\prod^\Gamma N_i/U$.  This is true because, by the $\eff$-elementarity of each $f_i$, 
$$M_i \vDash \neg \phi^j_k(m_i) \implies N_i \vDash \neg \phi^j_k(f_i(m_i))$$
So $k([i\mapsto m_i]_U)$ is a witness for $[i \mapsto f_i(m_i)]_U$.  Thus $f$ is a $K$-embedding.\hfill \dag\\

{\bf Proof of Theorem \ref{tameresult}:} We prove the type shortness and note that it implies the tameness by \cite[Theorem 3.5]{tamelc}.  Since we are not assuming amalgamation, we will show type shortness holds for atomic Galois equivalence.  Suppose that $X = \seq{x_i \in M_1 : i \in I}$ and $Y = \seq{y_i \in M_2:i \in I}$ are given such that, for all $I_0 \in P_\omega I$,

$$(\seq{x_i:i\in I_0}/\emptyset; M_1) E_{AT} (\seq{y_i:i\in I_0}/\emptyset; M_2)$$

That is, there is $N_{I_0} \in K$ and $f^\ell_{I_0}:M_\ell \to N_{I_0}$ such that $f_{I_0}^1(x_i) = f_{I_0}^2(y_i)$ for all $i \in I_0$.  Let $U$ be a fine ultrafilter on $P_\omega I$.  Then, following \cite{tamelc}, set 
\begin{itemize}

	\item $N = \prod^\Gamma_{I_0\in P_\omega I} N_{I_0}/U$;
	
	\item $f^\ell:M_\ell \to N$ is given by $f^\ell(m) = [I_0 \mapsto f^\ell_{I_0}(m)]_U$

\end{itemize}

N is well-defined by hypothesis and $f^\ell$ is a $K$-embedding by Lemma \ref{blank}.  For each $i \in I$, $\{I_0 \in P_\omega I: f^1_{I_0}(x_i) = f^2_{I_0}(y_i) \}$ contains $[i]:= \{I_0 \in P_\omega I: i \in I_0\} \in U$ by the fineness.  So $f^1(x_i) = f^2(y_i)$ for all $i \in I$.  Then
$$(X/\emptyset; M_1) E_{AT} (Y/\emptyset; M_2)$$\hfill\dag\\

Now, following \cite{bgcoheir}, we can define two notions of coheir.  There are two because the syntactic notion of type from $\mathcal{F}$ and the Galois notion of type from considering $K$ as an AEC do not necessarily coincide, although having the same Galois type implies having the same $\eff$-type

The first is Galois coheir $\dnfg$ (this could also be called $s$-coheir).  In this case, we consider Galois types over finite domains.  When Galois types are syntactic, these are complete syntactic types over a finite set.  The second is $t$-coheir $\dnft$, which is more like the first order version.

\begin{defin}
\begin{enumerate}
	\item Given $A, B, C \subset M$, we say $A \dnfg_C^{M} B$ iff
\begin{center}
for all finite $a \in A, b \in B, c \in C$, $gtp(a/bc)$ is realized in $C$.
\end{center}

	\item $K$ has the weak Galois order property iff there are finite tuples $\seq{a_i, b_i \in M : i < \omega}$ and $c$ and types $p \neq q \in gS(c)$ such that, for all $i, j < \omega$,
	$$j < i \implies a_i b_j \vDash q$$
	$$j \geq i \implies a_i b_j \vDash p$$

	\item Given $A, B, C \subset M$, we say $A \dnft_C^{M} B$ iff
\begin{center}
for all finite $a \in A, b \in B, c \in C$ and $\phi(x, y, z) \in \eff_K$, if $M \vDash \phi(a, b, c)$, then there is $c' \in C$ such that $M \vDash \phi(c', b, c)$.
\end{center}

	\item $K$ has the weak order property iff there are finite tuples $\seq{a_i, b_i \in M : i < \omega}$ and a formula $\phi(x, y, c) \in \eff_K$ with $c \in M$ such that, for all $i, j < \omega$,
	$$j < i \iff M \vDash \phi(a_i, b_j, c)$$
 
\end{enumerate}
\end{defin}

 Note that we have begun talking about Galois types over sets (rather than models, as standard) even though we only have amalgamation over models.  This adds some additional dificulties, but we are careful to avoid them here.  The adjective `weak' in describing the order property means that we only require $\omega$ length orders, rather than all ordinal lengths as in Shelah \cite{sh394}.

This ultraproduct allows us to weaken the requirements on getting this to be an independence relation the same way as in \cite[Section 8]{bgcoheir}.

\begin{theorem}\label{5.1g}
If $K$ is an averageable class with amalgamation that doesn't have the weak Galois order property and every model is $\aleph_0$-Galois saturated, then $\dnfg$ is an independence relation in the sense of \cite{bgcoheir}.
\end{theorem}

\begin{theorem}\label{5.1t}
If $K$ is an averageable class with amalgamation that doesn't have the weak order property and $\eff$ is first-order logic, then $\dnft$ is an independence relation in the sense of \cite{bgcoheir}.
\end{theorem}

Note that neither of these coheir's are precisely the definition given in \cite{bgcoheir}: as is standard, \cite{bgcoheir} only considered Galois types over models (so $\dnfg$ was not used) and there was no logic to choose (so $\dnft$ was not possible).  Nonetheless, the proofs of the above theorems go through the same arguments as in \cite[Theorem 5.1]{bgcoheir}.  The changes are minor, so we omit the details.  The interested reader can find the details on the author's website \cite{avcoheirdeets}; the above results are Theorems 6 and 10, respectively, from there.  Note that an advantage of using $\dnft$ is that Existence holds for free when $\eff$ is closed under existentials, although the disadvantadge is that $\dnft$ doesn't always have the semantic consequences often desired when dealing with types, i. e., if working in a class where syntactic types are not Galois types.  Additionally, with a little more stability, \cite{bgkv} shows that the two notions are the same if all models are $\aleph_0$-Galois saturated.

%Some additional facts about $EC$ classes:

%\begin{prop}
%Let $K = (EC(T, \Gamma), \prec)$ be an averageable\footnote{Not averageable class, but need that the type omitting part is a substructure.} class.
%\begin{itemize}
%	\item If $\prec$ is substructure and $T$ is a complete theory, then $K$ has amalgamation.
%	\item If $T$ is a complete theory with quantifier elimination, then $K$ has amalgamation and Galois types are syntactic.
%\end{itemize}
%\end{prop}

%{\bf Proof:} First, let $M_0 \prec M_1, M_2$ from $EC(T, \Gamma)$.\footnote{WB: Finish}\\ 

We have so far seen that averageable classes are very much like elementary classes.  The following result is a further restriction on the behavior of averageable classes.  It is easy to construct an averageable class with only a single model; take the standard model of arithmetic.  For general nonelementary classes, there are many more possibilities for the spectrum function of a class without arbitrarily large models.  However, the following result shows that, in the case of averageable classes, there are not.

\begin{theorem}\label{divline-thm}
Let $\Gamma$ be a finite set of countable existential\footnote{If $\Gamma$ consists of just quantifier free types, then the requirement in $(a)$ can be relaxed to just the same $\exists$- and $\forall$-theory.} types and let $M$ be a structure omitting $\Gamma$ that is $\Gamma$-closed.  Then, either
\begin{enumerate}
	\item[$(a)$] every $L$ structure omitting $\Gamma$ and satisfying the same $\exists\forall$-theory as $M$ is isomorphic to $M$; or
	\item[$(b)$] there are $\subset$-extensions of $M$ of all sizes, each satisfying the same $\exists\forall$-theory.
\end{enumerate}
\end{theorem}

We have stated the theorem in the simplest case.  However, variations are possible that strengthen the amount of \L o\'{s}' Theorem that holds and strengthen the similarity between the models; this means that it can be applied to situations such as $DLOG\mathbb{Z}$ or torsion modules over PIDs.  However, the countability remains crucial for the proof.

\begin{theorem} \label{divline2-thm}
Let $\Gamma$ be a finite set of countable types and let $M$ be a structure omitting $\Gamma$ that is $\Gamma$-nice.  Then, either
\begin{enumerate}
	\item[$(a)$] every $L$ structure omitting $\Gamma$ and elementarily equivalent to $M$ is isomorphic to $M$; or
	\item[$(b)$] there are $\prec$-extensions of $M$ of all sizes.
\end{enumerate}
\end{theorem}

{\bf Proof of Theorem \ref{divline-thm}:}  Enumerate each $p \in \Gamma$ as $\{\phi^p_n(x):n < \omega\}$.  Set $\psi_\ell(x) := \wedge_{p \in \Gamma} \vee_{n<\ell} \neg \phi^p_n(x)$.  We use these formulas to measure the type omission of all types of $\Gamma$ jointly.  Recall from Theorem \ref{propext}, that the $\Gamma$-ultraproduct produces a proper extension if there is an infinite subset of $M$ that all satisfy the same $\psi_\ell$.  This property separates our cases.

First, suppose this property fails; we will show that $(a)$ holds.  For each $\ell<\omega$, $\psi_\ell(M)$ is finite.  Thus, $M$ is countable and we can enumerate it as $\{m_i : i < \omega\}$.  For each $i < \omega$, pick some $\ell_i$ such that $M \vDash \psi_{\ell_i}(m_i)$.  Then, define
\begin{eqnarray*}
n_\ell &=& |\{n<\omega : \ell_i = \ell\}|\\
N_\ell &=& |\psi_\ell(M)|
\end{eqnarray*}
Note $n_\ell \leq N_\ell < \omega$.

Let $N$ be a model omitting $\Gamma$ and having the same $\exists\forall$-theory as $M$.  Note that
$$\exists x_0, \dots, x_{N_\ell-1} \left(  \wedge_{i < N_\ell} \psi_\ell(x_i)  \right)$$
$$\forall x_0, \dots, x_{N_\ell} \left(  \wedge_{i \leq N_\ell} \psi_\ell(x_i) \to \vee_{i \neq j \leq N_\ell} x_i = x_j  \right)$$
are both $\exists\forall$-sentences, so $|\psi_\ell(N)|=|\psi_\ell(M)|$.  We want to define bijections between these sets that fit together to be an isomorphism; this is done through a finite injury-style argument.

We construct sequences $\{n_i^L : i < L\}$ for $L < \omega$ such that
\begin{enumerate}
	\item for all $L < \omega$, $tp_{qf}(m_i : i < L) = tp_{qf}(n_i^L : i < L)$ and, for all $\ell < \omega$, 
	$$M \vDash \psi_\ell(m_i) \iff N \vDash \psi_\ell(n_i^L)$$
	
	\item for all $i < \omega$, $\seq{n_i^L : i < L < \omega}$ eventually stabilizes (in fact, it changes at most $\prod_{k \leq \ell} n_{\ell_k}$-many times)
\end{enumerate}

{\bf This is enough:}  For $i < \omega$, set $n_i$ to be the eventual value of $\seq{n_i^L : i < L < \omega}$.  Our isomorphism $f$ will take $m_i$ to $n_i$.  This is an isomorphism onto its range by the first part of $(1)$.  Furthermore, by the second part of $(1)$,
$$|\psi_\ell(M)| = |\psi_\ell(N) \cap f(M)| = |\psi_\ell(N)|$$
Since the $\psi_\ell(N)$ are finite and exhaust $N$, we have $f(M) = N$.  Thus, $M \cong N$, as desired.

{\bf Construction:} The following claim is key.

{\bf Claim:} For all $\bm \in M$, there is $\bn \in N$ such that $tp_{qf}(\bm) = tp_{qf}(\bn)$ and $N \vDash \psi_{\ell_i^*}(n_i)$, where $\ell_i^*$ is the picked witness for the $i$th member of $\bm$.

Suppose not.  Let 
$$N^* = \{ \bn' \in {}^{\ell(\bm)}N : \forall i . N \vDash \psi_{\ell_i^*}(n'_i) \}$$
Note that $N^*$ is finite.  Then, for each $\bn' \in N^*$, there is a quantifier-free $\phi_{\bn'}(\bx)$ that holds of $\bm$, but not of $\bn'$.  Set
$$\psi := ``\exists \bx \left( \bigwedge_i \psi_{\ell_i^*}(x_i) \wedge \bigwedge_{\bn' \in N^*} \phi_{\bn'}(\bx) \right)"$$
This is an $\exists(\neg \Gamma)$-sentence satisfied by $M$ and not by $N$, a contradiction.  Thus, the claim is proved.

Now we are ready to build $\{n_i^L : i < l\}$ by induction on $L<\omega$.

Set $n_0^1$ to satisfy the same qf-type as $m_0$ and satisfy the appropriate $\psi_\ell$.

For $L>1$, the above Claim says that there is at least one sequence satisfying $(1)$ for $m_0, \dots, m_{L-1}$.  Pick $\{n_i^L : i < L\}$ to be the sequence satisfying $(1)$ that agrees with the largest possible initial segment of $\{n_i^{L-1} : i < L-1\}$.

It is clear that this construction satisfies $(1)$.  To see it satisfies $(2)$, note that there are only finitely many choices for the $i$th element.  Thus, if an initial segment changed infinitely often, it would necessarily repeat; however, repetition is forbidden by the construction.\\

Second, suppose this property fails; we will show that $(b)$ holds.  We know that the $\Gamma$-ultraproduct is a proper extension.  We will iterate this.

For each ordinal, we will construct $M_\alpha \equiv_{\exists(\neg \Gamma)} M$ that omits $\Gamma$ and a coherent set of nonsurjective, embeddings $f_{\beta, \alpha} : M_\beta \to \alpha$ for $\beta < \alpha$.

For $\alpha = 0$, set $M_0 = M$.

For $\alpha = \beta+1$, set $M_\alpha := \prod^{\Gamma} M_\beta/U$, for $U$ a nonprincipal ultrafilter on $\omega$.  Note that this is a structure since the data is $\Gamma$-closed by Proposition \ref{gamclofo-prop} and Proposition \ref{comthe}.(2).  Then the ultrapower map $\j$ is a nonsurjective embedding.  Set $f_{\gamma, \alpha} = \j \circ f_{\gamma, \beta}$.

For $\alpha$ limit.  Let $U$ be a nonprinicipal uniform ultrafilter on $\alpha$ and set $M_\alpha := \prod^{\Gamma}_\beta M_\beta / U$; note that this is a $\Gamma$-ultraproduct rather than a $\Gamma$-ultrapower.  Again, this is a structure.  This shares the same $\exists (\neg \Gamma)$ theory of the $M_\beta$'s.

Define $f_{\beta, \alpha}:M_\beta \to M_\alpha$ by $f_{\beta, \alpha}(m) = [g^m_{\beta, \alpha}]_U$ where
$$g^m_{\beta, \alpha} (\gamma) = \begin{cases}
f_{\beta, \gamma}(m) & \beta \leq \gamma < \alpha \\
0 & \gamma < \beta
\end{cases}$$
Then this is a nonsurjective $K$-embedding such that $f_{\beta, \alpha} = f_{\gamma, \alpha} \circ f_{\beta, \gamma}$.

Since this chain is increasing, $M_\alpha \geq |M| + |\alpha|$, giving us the desired result. \hfill \dag

\section{Examples} \label{example-sec}

We now give several examples of classes $EC(T, \Gamma)$ for which our construction gives some compactness results.  The meaning of ``some compactness results'' is left vague, but the general behavior is that these are averageable classes for the appropriate fragment $\mathcal{F}$.  Another class of examples from torsion modules over PIDs is discussed in the next section.

As a final cautionary example, we discuss the case of Archimedean fields.  Typically, Archimedean fields are presented as ordered fields omitting the type of an infinite element $p_\infty(x) = \{ x > n \cdot 1 : n < \omega\}$.  However, if we take the theory of fields (of characteristic 0) and this type, then the data is not even $p$-closed: the $p_\infty$-ultraproduct has no positive infinite element, but does have infinitesimals and a negative infinite element; thus it's not closed under the field operations.  Thus, to fit into this framework, $\Gamma$ must contain continuum many types, one each to explicitly omit the positive and negative infinite elements and the infinitesimal elements above and below each standard element.  After these types are added, the class is $\Gamma$-closed with $\mathbb{R}$ as a maximal model of size $2^{\aleph_0} = |\mathbb{R}|$ (this maximality agrees nicely with Theorem \ref{propext}).

Another example along these lines is to consider differentially closed fields where every element is differentially algebraic over the constants (so it omits the type of a differential transcendental).

\subsection{Banach Spaces} \label{bs-ex} Banach spaces are the motivating example from this work: viewing continuous first-order logic as a certain fragment of $L_{\omega_1, \omega}$ (see Boney \cite{tdense}) lead to viewing the Banach space ultraproduct as one that, in part, omits unbounded elements by simply excluding them.  We outline how this can be put into this framework.

Let $L_b = \seq{B, R; +_B, 0_B; +_R, \cdot_R, 0_R, 1_R, <_R, c_r; \|\cdot\|, \cdot_{scalar}}_{r \in \mathbb{R}}$ be the two sorted language of normed linear spaces.  Then $T_b$ says that
\begin{itemize}
	\item $\{ c_r : r \in \mathbb{R}\}$ is a copy of $\mathbb{R}$; and
	\item $B$ is a vector space over $R$, with norm $\| \cdot \|:B \to R$.
\end{itemize}

We want to ensure that, in the ultraproduct, $R$ and $B$ each have no nonstandard elements, i. e., omit the type of an element of $R$ that is not some $c_r$.  Similar to the case of Archimedean fields, it is not enough to omit a single type; instead every nonnegative real must have a types specifying there is no nonstandard real around it and a type specifying there are not Banach space elements that would be mapped to such an element.

\begin{itemize}
	\item $p_{\infty}(x) = \{ R(x) \wedge (x < -n \vee n < x) : n < \omega\}$;
	\item $p_r(x) = \{R(x) \wedge (x \neq c_r) \wedge (c_{r -\frac{1}{n}} < x< c_{r+\frac{1}{n}}) : n < \omega\}$ for $r \in \mathbb{R}$;
	\item $q_{\infty}(x) = \{ B(x) \wedge (\|x\| < -n \vee n < \|x\| ) : n < \omega\}$; and
	\item $q_r(x) = \{ B(x) \wedge (\|x\| \neq c_r) \wedge (c_{r - \frac{1}{n}} < x < c_{r+\frac{1}{n}}) : n < \omega\}$.
\end{itemize}

Set $\Gamma = \{ p_r (x) : r \in \mathbb{R} \cup \{ \infty\}\} \cup \{ q_r(x) : r \in \mathbb{R}^{\geq 0} \cup \{\infty\}\}$.  We omit the details, but $EC(T_b, \Gamma)$ is $\Gamma$-closed: the key details is that the standard real number that two sequences correspond to can be used to calculate the standard real number their sum or product corresponds to.  This means that the Universal \L o\'{s}' Theorem holds.  Additionally, by Observation \ref{doublestar}, the class of formulas which \L o\'{s}' Theorem holds is closed under ``bounded quantification,'' that is, of the form
$$\exists x \left( \phi(x, \by) \wedge \|x\| < c \right)$$
for some $c > 0$.

Comparing this with first-order continuous logic, there is not a requirement that the space be of bounded diameter.  Moreover, the condition above recovers some of the results from Ben Yaacov \cite{byunbounded} about unbounded metric spaces.

Other results for continuous logic can be recovered through these methods.  For instance, when trying to extend the language $L_b$ and preserve the $\Gamma$-closedness of the class, the relevant condition turns out to be uniform continuity of the function or relation, which agrees with the results from continuous first-order logic.  Additionally, a continuous version of the $\Gamma$-ultraproduct can be developed along the same lines.

%$\Gamma(\mathcal{B})$ is an $L_b$-structure for any $\mathcal{B}$.  We omit the details, but standard real number that two sequences correspond to, which are their witnesses to inclusion in $\Gamma(\mathcal{B})$, can be used to calculate which number their sum or product corresponds to.  

%Considering the class $EC(T, \Gamma)$, the models are normed vector spaces and, following standard proofs, $\Gamma\left(\prod \mathcal{B}_i/U\right)$ will be complete and satisfy the Universal \L o\'{s}' Theorem.  Moreover, if $\phi(x, \by)$ is a formula that transfer from a $U$-large set to $\Gamma\left(\prod^\Gamma \mathcal{B}_i/U\right)$, then so is
% $$\exists x \left( \phi(x; \by) \wedge \|x\| < c_r\right)$$
%for each $r\in \mathbb{R}$ be Observation \ref{doublestar}\footnote{Track changes\\WB: Add this observation back!}.  This recovers some results from \cite{byunbounded}.

%Note, that in contrast to continuous first-order logic (Ben Yaacov, Berentstein, Henson, Usvyatsov \cite{fourguys} is the standard reference), we do not require that the 

%the metric spaces under consideration have bounded diameter, so this extra step is not necessary.  However, Ben Yaacov has explored the context of infinite diameter metric spaces in Ben Yaacov \cite{byunbounded} and there are many similarities, especially our criteria for existential closure in Observation \ref{doublestar} and the criteria for closure under $\inf$ in his logic

%EXTENDING THE LANGUAGE AND UNIFORM CONTINUITY (CAICEDO, EAGLE, IOVINO)

\subsection{Abelian Torsion Groups} \label{atg-ex} Let $L_g = \{+, 0, -\cdot\}$ and $T_{ag}$ be the theory of abelian groups.  Abelian torsion groups are models of $T_{ag}$ that omit $tor(x) = \{ n \cdot x \neq 0 : n < \omega\}$.  We claim that abelian torsion groups are $tor$-closed.

\begin{prop}
If $G$ is an abelian group, then it is $tor$-closed.
\end{prop}

{\bf Proof:} Given $g \in G$, we have that $G \vDash \neg (n \cdot g \neq 0)$ exactly when $o(g) \mid n$.  Since $o(g) = o(-g)$ and $o(g_1 + g_2) = \lcm(o(g_1), o(g_2)) \mid o(g_1) o(g_2)$, setting $g_-(n) = n$ and $g_{+}(n, m) = nm$ shows that $G$ is $tor$-closed. \hfill \dag\\

A more in depth analysis shows the full \L o\'{s}' Theorem holds in the wider class of torsion modules over a PID.

\subsection{DLOG$\mathbb{Z}$} \label{dlogz-ex} We consider the theory of densely ordered abelian groups\footnote{Note that the group structure is not crucial here, and the same analysis could be done with the theory of dense linear orders with a cofinal $\mathbb{Z}$-chain.} with the infinitary property of having a cofinal $\mathbb{Z}$-chain.  The first order part of this theory was first shown to have quantifier elimination by Skolem \cite{skolem}\footnote{For a little more history, see the introduction of Hieronymi \cite{hiero}.  Also, Miller \cite{miller} contains a proof and is more easily accessible than Skolem's original}.  We will show that the first order portion of the theory is preserved by the appropriate $p$-ultraproduct, and then use quantifier elimination to bootstrap the full version of \L o\'{s}' Theorem.

Set $T := Th(\mathbb{Q}, <, +, -, 0, 1, n)_{n  \in \mathbb{Z}}$ and $Z(x) := \{ x \leq c_n \text{ or } c_m \leq x : n < m \in \mathbb{Z}\}$, where $c_n$ is the constant representing $n$.  By a model of DLOG$\mathbb{Z}$, we mean a model of $T$ that omits $Z$, i. e. one where $\{c_n : n \in \mathbb{Z}\}$ is a discrete, countable sequence that is cofinal in both directions. This theory has quantifier elimination and is axiomatized by the axioms for an ordered, uniquely divisible, torsion-free abelian group that is dense as an ordering and the elementary diagram of $(\mathbb{Z}, +, <)$.

\begin{prop}
$(EC(T, Z), \prec)$ is closed under $Z$-ultraproducts and they satisfy \L o\'{s}' Theorem.  Morever, this is a class with amalgamation where Galois types are syntactic.
\end{prop}

{\bf Proof:}  Let $M_i$ be a model of DLOG$\mathbb{Z}$ for each $i \in I$ and let $U$ be an ultrafilter on $I$.

\begin{claim}
$\prod^Z M_i/U$ is a structure that models $T$.
\end{claim}

{\bf Proof:} We have to show that it contains the constants and is closed under functions.  Each $c_n$ is represented by $[i \mapsto c_n^{M_i}]_U$, which fails to satisfy ``$x \leq c_{n-1} \text{ or } c_{n+1} \leq x$'' everywhere.  Next we look at addition; subtraction is similar.  Let $[f]_U, [g]_U \in \prod^Z M_i/U$ that are witnessed by 
\begin{eqnarray*}
c_{n_f} < &[f]_U& < c_{m_f}\\
c_{n_g} < &[g]_U& < c_{m_g}
\end{eqnarray*}
Then $[f]_U + [g]_U = [f+g]_U \in \prod^Z M_i/U$ as witnessed by
\begin{eqnarray*}
c_{n_{f+g}} = c_{n_f} + c_{n_g} < [f]_U + [g]_U < c_{m_f} + c_{m_g} = c_{m_{f+g}}
\end{eqnarray*}
Since $\prod^Z M_i/U$ is a structure, we now wish to show it models $T$.  We know $\exists \forall$-sentences transfer, so we only need to show that the existentials in the divisibility of the group and denseness of the order hold.\\
For the divisibility, suppose $[f]_U \in \prod^Z M_i/U$ and $k < \omega$ such that there is $X \in U$ and $n_f < m_f \in \mathbb{Z}$ such that, for all $i \in X$, $M_i \vDash c_{n_f} < f(i) < c_{n_g}$.  Then, for each $i \in I$, there is $\frac{f}{k}(i) \in M_i$ such that 
$$M_i \vDash k \cdot \frac{f}{k}(i) = f(i) \wedge (c_{-|n_f| - |m_f|} < \frac{f}{k}(i) < c_{|n_f| + |m_f|}$$
For the denseness, suppose $[f]_U, [g]_U \in \prod^Z M_i/U$ such that $\prod^Z M_i/U \vDash [f]_u < [g]_U$.  Thus, there is $X \in U$ and $n_f <  m_f, n_g<m_g \in \mathbb{Z}$ such that, for all $i \in X$, we have
\begin{enumerate}
	\item $M_i \vDash f(i) < g(i)$;
	\item $M_i \vDash c_{n_f} < f(i) < c_{m_f}$; and
	\item $M_i \vDash c_{n_g} < g(i) < c_{m_g}$.
\end{enumerate}
We can find $h \in \prod M_i$ such that $M_i \vDash f(i) <h(i) < g(i)$.  For $i \in X$, we have $M_i \vDash c_{n_f} < h(i) < c_{n_g}$.  Thus, $[h]_U \in \prod^Z M_i/U$ and $\prod^Z M_i / U \vDash [f]_U < [h]_U < [g]_U$.

\begin{claim}
The $Z$-ultraproduct satisfies \L o\'{s}' Theorem.
\end{claim}

%This easily follows from the above\footnote{WB: This proof should probably exist above as a general statement}.  Let $\phi(\bx)$ be a formula and $[f_0]_U, \dots, [f_{n-1}]_U \in \prod^p M_i/U$.  By quantifier elimination, there is a quantifier free $\psi(\bx)$ such that $\phi$ and $\psi$ are equivalent modulo $T$.  Then 
%\begin{eqnarray*}
%\{i \in I : M_i \vDash \phi(f_0(i), \dots, f_{n-1}(i)) \} \in U &\iff& \{i \in I : M_i \vDash \psi(f_0(i), \dots, f_{n-1}(i)) \} \in U\\
%&\iff& \prod^p M_i/U \vDash \psi([f_0]_U, \dots, [f_{n-1}]_U) \\
%&\iff& \prod^p M_i/U \vDash \phi([f_0]_U, \dots, [f_{n-1}]_U)
%\end{eqnarray*}
This follows by Proposition \ref{qelos-prop} and quantifier elimination.

Second, we show that the class has amalgamation and that Galois types are syntactic.  Note that, by definition of the class, having the same Galois type implies having the same syntactic type.  Let $M_0 \prec M_1, M_2 \in EC(T, Z)$, possibly with $a_\ell \in M_\ell$ such that $tp(a_1/M_0; M_1) = tp(a_2/M_0; M_2)$.  Then, since the elementary class of models of $T$ has amalgamation and has that syntactic types are Galois types, there is $N^* \vDash T$ and $f_\ell: M_\ell \to_{M_0} N^*$ such that, if we are dealing with types, $f_1(a_1) = f_2(a_2)$.  $N^*$ might realize $p$, but set $N$ to be the substructure of $N^*$ with universe $\{ x \in N  \mid \exists n, m \in \mathbb{Z}. N^* \vDash c_n < x < c_m\}$.  This is a substructure of $N^*$ that models $T$, contains $f_1(M_1)$ and $f_2(M_2)$, and omits $p$.  By quantifier elimination, these inclusions are actually elementary substructure.  Thus $N$ is the desired amalgam.  Additionally, if we are dealing with the type statement, $f_1(a_1) = f_2(a_2)$, so $gtp(a_1/M_0; M_1) = gtp(a_2/M_0; M_2)$ as desired. \hfill \dag

This example can be generalized by looking at ordered $R$-vector spaces over an ordered division ring $R$ rather than just ordered divisible abelian group.  By \cite[Corollary 1.(7.8)]{ominbook}, this wider class also has quantifier elimination and the argument works in the same way.

\subsection{Multi-sorted first order logic} \label{sort-ex}

Take a multi-sorted language $\tau$ with sorts $\{S_\alpha : \alpha < \kappa\}$ and a theory $T$.  There is a natural correspondence between multi-sorted models of $T$ and models of a (non-sorted) first-order theory $T^*$ in the language $\tau^* := \tau \cup \{S_\alpha : \alpha < \kappa\}$ that omit the type $sort(x) := \{ \neg S_\alpha (x) : \alpha < \kappa\}$.  Then, the class $EC(T^*, sort)$ is not $sort$-nice, but is still well-behaved with respect to the $sort$-ultraproduct in the following sense.

\begin{prop}
$EC(T^*, sort)$ satisfies \L o\'{s}' Theorem with respect to $\tau^*$ formulas that come from sorted $\tau$ formulas.
\end{prop}

{\bf Proof:}  First, we observe that the class is $sort$-closed: if $F$ be a function of $\tau$, then $T^*$ determines the sort of $F$ applied to any valid input.  This means that a the universal \L o\'{s}' Theorem holds.  Moreover, suppose that $\exists x \phi(x, \by)$ is a $\tau^*$ formula that comes from a sorted $\tau$ formula.  Then, this formula determines which sort a witness $x$ would be in.  This is precisely the information required to define the function $g_{\exists x \phi(x, \by)}$; note that it is a constant function.  Thus, the set of formulas that $EC(T^*, sort)$ satisfies \L o\'{s}' Theorem with contain all quantifier-free $\tau^*$ formulas that come from sorted $\tau$ formulas and is closed under exstentials.  By applying Remark \ref{inductlos-remark}, this extends to the class of all formulas coming from sorted $\tau$ formulas, as desired. \hfill \dag

This allows one to read off the normal compactness results of sorted first-order logic from the results of this paper; moreover, Theorem \ref{divline-thm} means that a sorted model has proper elementary extensions iff one of the sorts has infinite size.  Indeed, this correspondence seems to go both ways and one could likely perform the same analysis in this paper by looking at which $EC$ classes are ``sortable.''

\subsection{Highly Complete Ultrafilters} Our final example shows that, if there are very complete ultraproducts, then this new ultraproduct coincides with the classic one.

\begin{theorem}
If $U$ is $\chi$-complete, $\chi > |\Gamma|$, and $\chi > |p|$ for all $p \in \Gamma$, then $\prod^\Gamma M_i/U = \prod M_i/U$.
\end{theorem}

{\bf Proof:} We always have $\prod^\Gamma M_i \subset \prod M_i$.  Let $f \in \prod M_i$.  We want to show $f \in \prod^\Gamma M_i$ by finding a witness.  For each $\phi \in p \in \Gamma$, set $X_\phi^{f, p}:= \{ i \in I : M_i \vDash \neg \phi(f(i))\}$.  For each $p \in \Gamma$, $I$ is the union of $\{ X_\phi^{f, p} : \phi \in p\}$.  Since $|p| < \chi$, there is some $\phi_p$ such that $X_{\phi_p}^{f,p} \in U$. Then
$$X^f = \cap_{p \in \Gamma} X^{f, p}_{\phi_p} \in U$$
shows that the map $p \mapsto \phi_p$ is a witness.  Thus $\prod^\Gamma M_i = \prod M_i$. \hfill \dag\\

Note that, if $\kappa$ is some large cardinal giving rise to $\kappa$-complete ultrafilters and $K$ is averageable with respect to $\kappa$-complete ultrafilters, then $K$ will satisfy the relevant parts of the last section with $\kappa$ in place of $\omega$; see \cite{tamelc} and \cite[Section 8]{bgcoheir} for what is relevant.

\section{Torsion Modules} \label{tormod-sec}

In this section, we explore the previous results applied to torsion modules over PIDs and apply some results for nonelementary stability theory.  The stability theoretic results are not deep (and probably follow from results about modules and other properties of torsion modules), but we intend this to show what can be done.

\subsection{The Torsion Ultraproduct}

For this subsection, assume that $R$ is a commutative ring with unity.\footnote{The following weakening of commutativity is also sufficient: $\forall x \forall y \exists z (xy = zx)$.  Then we can take the ultraproduct of left torsion modules.}

We review some basics of the model theory of modules, using Prest \cite{prestmodulebook} as the reference.  The language is $L_R = \seq{+, r\cdot, -, 0}_{r \in R}$.  Then the theory of $R$-modules $T_R$ is the statement of all of the module axioms; note that this is a universal theory.  Given a module $M$ and $m \in M$, set
$$\mathcal{O}^M(m) := \{ r \in R : r \cdot m = 0 \text{ and } r \text{ is regular}\}$$
Recall that regular elements are those that are not zero divisors.  We drop the $M$ if it is clear.  If this set is non-empty and $m \neq 0$, then $m$ is a torsion element and every element of $\mathcal{O}(m)$ is called an order of $m$.  If every element of $M$ is a torsion element, then $M$ is a torsion module.

Note that Shelah \cite{shinfmod} has recently explored the more general behavior of $L_{\lambda, \mu}$-theories of modules, but does not deal with compactness or nonforking\footnote{\cite{shinfmod} says he intends to deal with nonforking in \cite{shF1210}, but this has yet to appear.}.

Set $tor(x) = \{ r \cdot x \neq 0 : r \in R\}$ to be the type of a torsion-free element.  Let $\{M_i : i \in I\}$ be a collection of torsion modules (i.e. modules that omit $tor$) and let $U$ be an ultrafilter on $I$.  Then the $tor$-ultraproduct is
\begin{eqnarray*}
\prod^{tor} M_i /U := \{ [f]_U &:& f \in \prod M_i \te{ and there is } X_f \in U \te{ and } r_f \in R\\ & & \te{ such that } r_f \in \mathcal{O}^{M_i}(m_i) \te{ for all } i \in X_f \}
\end{eqnarray*}

\begin{prop}  \label{struct-prop}\label{struct-thm}
$EC(T_R, tor)$ is $tor$-closed and the universal \L o\'{s}' Theorem holds.
\end{prop}

{\bf Proof:} Note that the first part implies the second by Theorem \ref{weaklos}.

We need to construct functions that tell us the order of a sum, etc. based on the order of the inputs.  For later use, we do more: for each $L_R$-term $\tau(\bx)$, we inductively construct $f_\tau:R^{\ell(\bx)} \to R$ such that $f_\tau''\prod\mathcal{O}(m_i) \subset \mathcal{O}(\tau(m_0, \dots, m_{n-1}))$
\begin{itemize}
	\item if $\tau(\bx) = x_i$, then $f_\tau(\br) = r_i$;
	\item if $\tau = s \cdot \sigma$, then $f_\tau = f_\sigma$;
	\item if $\tau = \sigma + \chi$, then $f_\tau = f_\sigma f_\chi$; and
	\item if $\tau = - \sigma$, then $f_\tau = f_\sigma$.\footnotei{WB: Do we use this later?\\Gets used at least in Lemma \ref{pplos-lem}}
\end{itemize}
Thus, $EC(T_R, tor)$ is $tor$-closed..

%\begin{prop} \label{struct-prop}\label{struct-thm}
%Let $\{M_i : i \in I\}$ be a collection of torsion modules and $U$ is an ultrafilter on $I$.  Then $\prod^{tor} M_i / U$ is a torsion module and the following universal variant of \L o\'{s}' Theorem holds: if $[f_0]_U, \dots, [f_{n-1}]_U \in \prod^{tor}M_i/U$ and $\phi(\bx)$ is a universal formula in $L_R$, then
%$$\{ i \in I : M_i \vDash \phi(f_0(i), \dots, f_{n-1}(i)) \} \in U \implies \prod^{tor} M_i / U \vDash \phi([f_0]_U, \dots, [f_{n-1}]_U)$$
%\end{prop}

%{\bf Proof:}  First, we show that $\prod^{tor} M_i / U$ is an $L_R$-structure.  To do so, given $[f]_U, [g]_U \in \prod^{{tor}} M_i / U$ with orders $r_f, r_g \in R$, we need to compute orders of the terms of $R$ of $[f]_U$ and $[g]_U$.  For later use, we do more: for each $L_R$-term $\tau(\bx)$, we inductively construct $f_\tau:R^{\ell(\bx)} \to R$ such that $f_\tau''\prod\mathcal{O}(m_i) \subset \mathcal{O}(\tau(m_0, \dots, m_{n-1}))$
%\begin{itemize}
%	\item if $\tau(\bx) = x_i$, then $f_\tau(\br) = r_i$;
%	\item if $\tau = s \cdot \sigma$, then $f_\tau = f_\sigma$;
%	\item if $\tau = \sigma + \chi$, then $f_\tau = f_\sigma f_\chi$; and
%	\item if $\tau = - \sigma$, then $f_\tau = f_\sigma$.\footnotei{WB: Do we use this later?}
%\end{itemize}
%Thus, $\prod^{tor} M_i / U$ is an $L_R$-structure.  The universal variant of \L o\'{s}' Theorem then follows exactly from \cite[Theorem 2.3]{nkou}.  Since $T_R$ is universal and $tor$ is quantifier-free, this means that $\prod^{tor} M_i / U$ is a torsion module, see \cite[Proposition 2.7]{nkou}. \hfill \dag \\

In fact, in this case, we have $\prod^{tor} M_i/U$ is precisely the torsion subgroup of the full ultraproduct $\prod M_i/U$.  Thus, the construction of the $tor$-ultraproduct is not new, but we can use the results from earlier sections and the model theory of modules to get some new results.\footnotei{WB: Obvi, the goal is to remove the assumption of PIDs or find an example that shows it is necessary}

\footnotei{SAME $\exists \forall$ THEORY!!!!}

Further study of the ultraproduct requires specialization to PIDs, but we already have the following dividing line for modules.  Roughly, this says that, given a torsion module over a countable commutative ring, either it is the only torsion module like it \emph{or} there are torsion modules like it of all sizes.

\begin{cor} \label{dividingline-thm}
Let $M$ be a torsion module over a countable, commutative ring.  Then either
\begin{enumerate}
	\item every torsion module that is $\exists$- and $\forall$-equivalent to $M$ is in fact isomorphic to $M$; or
	\item there are torsion $\subset_{\forall}$-extensions of $M$ of all sizes (in fact, the all model the same $\exists\forall$ theory).
\end{enumerate}
\end{cor}

Note that there is no explicitly stated restriction on the size of the module in (1), but $M$ will necessarily be countable as will any torsion module $\forall$-equivalent to it.

{\bf Proof:} This is Theorem \ref{divline-thm} in this context.\hfill\dag\\

\subsection{Torsion Compactness over PIDs} \label{torcomp-sec}

For the remained of this subsection, assume that $R$ is a principal ideal domain.  Note that PIDs are integral domains, so all nonzero elements are regular.

The goal of this subsection is to prove \L o\'{s}' Theorem for elementarily equivalent modules.  The proof of this uses Proposition \ref{qelos-prop} and has two steps:

\begin{enumerate}
	\item recall that $T_R$ has p. p. elimination of quantifiers; and
	\item show that \L o\'{s}' Theorem holds for p. p. formulas (and a little more).
\end{enumerate}

We need to recall the key facts about p. p. elimination of quantifiers.

\begin{defin} \label{ppform-def}
$\phi(\bx)$ is a p. p. (primitive positive) formula iff it is a conjunction of formulas of the form $p^n \divides \tau(\bx)$ and $\tau(\bx)= 0$ for a term $\tau$, a prime $p \in R$, and $n < \omega$.
\end{defin}

Note that p. p. formulas have a more general definition (see \cite[Section 2]{prestmodulebook} for the more general definition and a deeper discussion of their role in the model thoery of modules), but this is an equivalent formulation in PIDs (\cite[Theorem 2.$\mathbb{Z}$1]{prestmodulebook}).  Indeed this formulation is the key reason we have specified to PIDs as it allows us to prove \L o\'{s}' Theorem for p. p. formulas.  Note that Shelah \cite[Theorem 2.4]{shinfmod} has a much more general version of this result for $L_{\lambda, \theta}$-theories of modules (note he calls these formulas p. e. or ``positive existential''), but the first-order version suffices.

Given a complete theory of a module, all formulas are equivalent to a boolean combination of p. p. formulas.  However, a more precise result involving invariants conditions is true.

\begin{defin}\label{invcond-def}
Given a module and p. p. formulas $\phi(\bx)$ and $\psi(\bx)$, set $Inv(M, \phi, \psi) = |\phi(M)/\phi(M) \cap \psi(M)|$.  An \emph{invariants condition} is the assertion that $Inv(M, \phi, \psi)$ is either greater than or less than some $k<\omega$.
\end{defin}

\begin{fact}[\cite{prestmodulebook}.2.13]\label{ppeq-fact}
If $\phi(\bx)$ is a formula, then there is a boolean combination of invariants conditions $\sigma$ and a boolean combination of p. p. formulas $\psi(\bx)$ such that 
$$T_R \vdash \forall \bx (\phi(\bx) \iff (\sigma \wedge \psi(\bx)))$$
\end{fact}

\begin{lemma} \label{pplos-lem}
Suppose all $M_i$ are elementarily equivalent.  Given $[f_0]_U, \dots, [f_{n-1}]_U \in \prod^{tor} M_i /U$ and p. p. $\phi(\bx)$,
$$\{ i \in I : M_i \vDash \phi(f_0(i), f_{n-1}(i)) \} \in U \iff \prod^{tor} M_i / U \vDash \phi([f_0]_U, \dots, [f_{n-1}]_U)$$
\end{lemma}

{\bf Proof:} Note that $\neg \phi(\bx)$ is universal, so right to left follows from Theorem \ref{struct-thm} above.  For the other direction, suppose $\phi(\bx)$ is of the following form:
$$\bigwedge_{j<m} (\exists y_j. p_j^{n_j} \cdot y_j = \tau_j(\bx) ) \wedge \bigwedge_{j<m'} (\sigma_j(\bx)=0)$$
and that $Y:=\{i \in I : M_i \vDash \phi(f_0(i), \dots, f_{n-1}(i)) \} \in U$.  The difficulty is establishing that the existential witnesses lie in the $tor$-ultraproduct.  By the definition of the $tor$-ultraproduct, each parameter $[f_k]_U$ has some fixed order on a $U$-large set, say $r_k \in \mathcal{O}^{M_i}(f_k(i))$ for all $i \in X_k \in U$.  Then $r:=f_{\tau_i}(r_0, \dots, r_{n-1})$ will be an order for them on $X:= \cap_{k < n}X_k$; recall that $f_\tau$ was constructed in the proof of Proposition \ref{struct-prop}.

For each $i \in Y$ and $j < m$, find $m^i_j$ such that 
$$M_i \vDash p_j^{n_j} \cdot m_j^i = \tau_j(f_0(i), \dots, f_{n-1}(i))$$
Then $r$ is also an order for each $m_j^i$ when $i \in X$.  We define $g_j \in \prod M_i$ by
$$g_j(i) = 
\begin{cases} 
	m_j^i &\text{if } i \in X \cap Y \\ 
	0 & \text{otherwise} 
\end{cases} 
$$
Then $r$ and $X \cap Y \in U$ witness that $[g_j]_U \in \prod^{tor} M_i/U$ for each $j < m$ and 
$$\prod^{tor} M_i/U \vDash p_j^{n_j} \cdot [g_j]_U = \tau_j([f_0]_U, \dots, [f_{n-1}]_U)$$
Thus, $\prod^{tor} M_i/U \vDash \phi([f_0]_U, \dots, [f_{n-1}]_U)$, as desired. \hfill\\

We can easily extend this result to boolean combinations of p. p. formulas.

\begin{cor} \label{boolpp-cor}
Suppose $M_i$ are elementarily equivalent.  Given $[f_0]_U, \dots, [f_{n-1}]_U \in \prod^{tor} M_i/U$ and a boolean combination of p. p. $\phi(\bx)$,
$$\{i \in I : M_i \vDash \phi(f_0(i), \dots, f_{n-1}(i)) \} \in U \iff \prod^{tor} M_i/U \vDash \phi([f_0]_U, \dots, [f_{n-1}]_U)$$
\end{cor}

\begin{remark}
Lemma \ref{pplos-lem} is the key result that requires the specialization to modules over PIDs, and it's not currently known if this holds in general for commutative rings.  As an alternate hypothesis, this result also holds if all annihilator ideals are prime.
\end{remark}

We can use the fact that \L o\'{s}' Theorem holds for p. p. formulas to show that it also holds for boolean combinations of invariants conditions (these are sometimes called invariants sentences).

\begin{lemma}\label{invsentlos-lem}
Suppose that all $M_i$ are elementarily equivalent.  Let $\phi$ be a boolean combination of invariants conditions.  Then $\phi$ is part of the common theory of the $M_i$'s iff $\prod^{tor} M_i/U \vDash \phi$.
\end{lemma}

{\bf Proof:} Since a negation of a boolean combination is itself a boolean combination, it suffices to show one direction.  Thus, assume that $M_i \vDash \phi$ for all $i \in I$.  Since conjunctions and disjunctions transfer (see Remark \ref{inductlos-remark} or easy to work out the details), we only have to show this for $Inv(M, \phi, \psi) \geq k$ and $Inv(M, \phi, \psi) < k$.  WLOG, assume $\psi$ implies $\phi$.  Note that
\begin{eqnarray*}
Inv(M, \phi, \psi) \geq k &\equiv& ``\exists \bv_0, \dots, \bv_{k-1} ( \bigwedge_{i<k}\phi(\bv_i) \wedge \bigwedge_{j < i < k} \neg \psi(\bv_j - \bv_i))''\\
Inv(M, \phi, \psi) < k &\equiv& ``\forall \bv_0, \dots, \bv_{k-1} (\bigvee_{i < k} \neg\phi(\bv_i) \vee \bigvee_{j<i<k} \psi(\bv_j - \bv_i))''
\end{eqnarray*}
$Inv(M, \phi, \psi) \geq k$ is $\exists\forall$, so the result holds by Proposition \ref{comthe}.  The formula ``$\bigvee_{i < k} \neg\phi(\bv_i) \vee \bigvee_{j<i<k} \psi(\bv_j - \bv_i)$'' is a boolean combination of p. p. formulas, so it transfers by Corollary \ref{boolpp-cor}.  Then $Inv(M, \phi, \psi) < k$ is universal over a formula that transfers, so it transfers as well; again, see Remark \ref{inductlos-remark} or work out the details. \hfill \dag\\

We now have all of the tools that we need to prove the full version of \L o\'{s}' Theorem.

\begin{theorem}\label{fulllos-thm}
Let $T_R^*$ be a complete theory $L_R$-theory extending $T_R$ (recall $R$ is a PID).  Then $EC(T_R^*, tor)$ satisfies \L o\'{s}' Theorem with the $tor$-ultraproduct.
\end{theorem}

{\bf Proof:}  Suppose that $\{M_i \mid i \in I\}$ are torsion models of $T_R^*$ and $U$ is an ultrafilter on $I$.  Then $\prod^{tor} M_i/U$ is a torsion module by Proposition \ref{struct-prop}. Let $[f_0]_U, \dots, [f_{n-1}]_U \in \prod^{tor}M_i/U$ and $\phi(\bx)$ be a formula.  By Fact \ref{ppeq-fact}, $\phi(\bx)$ is equivalent modulo $T_R$ to a boolean combination $\sigma$ of invariants conditions and a boolean combination $\psi(\bx)$ of p. p. formauls.  Then
\begin{eqnarray*}
\{ i\in I : M_i \vDash \phi(f_0(i), \dots, f_{n-1}(i)) \} \in U &\iff& \{i \in I : M_i \vDash \sigma \wedge \psi(f_0(i), \dots, f_{n-1}(i)) \} \in U\\
&\iff& \prod^{tor} M_i /U \vDash \sigma \wedge \psi([f_0]_U, \dots, [f_{n-1}]_U)\\
&\iff& \prod^{tor} M_i/U \vDash \phi([f_0]_U, \dots, [f_{n-1}]_U)
\end{eqnarray*}
The first and third equivalence is by Fact \ref{ppeq-fact} and the second equivalence is by Corollary \ref{boolpp-cor} and Lemma \ref{invsentlos-lem}.\hfill \dag\\

\subsection{Examples}

For this subsection, we specialize to $R = \mathbb{Z}$.  That is, we examine abelian torsion groups.

We look at some examples of torsion abelian groups and examine how the groups differ from their $tor$-ultraproducts and how the AEC $(\Mod(T \cup \{\forall x \bigvee_{n<\omega} n \cdot x = 0\}), \prec)$ differs from the elementary class $(\Mod(T), \prec)$.

We list some torsion abelian groups $G$ such that $G \precneqq \prod^{tor} G/U \precneqq \prod G / U$\footnote{Formally, $G$ is not a subset of $\prod^{tor}G/U$, but is canonically embedded in it; we blur this distinction by identifying $g \in G$ and $[i \mapsto g]_U\in \prod^{tor}G/U$}.  The main point here is the inequalities, as the elementary substructure results follow from Theorem \ref{fulllos-thm}.  If $G$ does not have finite exponent, then $\prod G/U$ necessarily contains elements with no order, so $\prod^{tor} G/U \subsetneq \prod G/U$.  Theorem \ref{propext} gives a  condition for $G \subsetneq \prod^{tor} G/U$.  In this context, the result becomes:
\begin{center}
If there is some $n < \omega$ such that there are infinitely many $g \in G$ such that $o(g) \divides n$, then $G \subsetneq \prod^{tor} G/U$.
\end{center}

Thus, the following groups are all proper elementary subgroups of their $tor$-ultraproducts (for any nonprincipal ultrafilter).
\begin{enumerate}
	\item $\oplus_{n<\omega} \mathbb{Z}_n$
	\item More generally, $\oplus_{n<\omega} \mathbb{Z}_{s_n}$ for any sequence $\seq{s_n : n < \omega}$ such that there is a prime $p$ that divides infinitely many of the $s_n$'s
	\item $\oplus_{n<\omega} \mathbb{Z}(p^\infty)$ for any prime $p$ or $\oplus_{n<\omega} \mathbb{Q}/\mathbb{Z}$
\end{enumerate}

Note that $\mathbb{Z}(p^\infty)$ and $\mathbb{Q}/\mathbb{Z}$ (or the sum of finitely many of them) do not satisfy this criterion.  Thus, by Theorem \ref{dividingline-thm}, any torsion group elementarily equivalent to them is in fact isomorphic to them.

Consider $G = \oplus_{n<\omega} \mathbb{Z}_{2^n}$ and $I = \omega$.  Note that, for any nonzero $g \in G$, there is a maximum $k < \omega$ such that $2^k \divides g$.  However, this is not the case in $\prod^{tor} G/U$: set $f:I \to \omega$ sucht that $f(i)$ is $2^{i-1}$ in $\mathbb{Z}_{2^i}$, i.e. $f(i) \in \prod \mathbb{Z}_{2^n}$ such that
$$f(i)(n) = 
\begin{cases} 
	2^{i-1} &\text{if } i = n \\ 
	0 & \text{otherwise} 
\end{cases} 
$$
Then each $f(i)$ has order 2, so $[f]_U \in \prod^{tor}G/U$.  However, for every $k < \omega$,
$$\{i \in I : G \vDash \exists y. 2^k \cdot y =f(i)\} = \omega - (k+1) \in U$$
So $\prod^{tor} G/U \vDash 2^k \divides [f]_U$ for all $k < \omega$.  Thus there are countable submodels of $\prod^{tor} G/U$ not isomorphic to $G$.  Thus, $Th(G)$ is not countably categorical amongst abelian torsion groups.

In contrast, we now examine $Th(\oplus \mathbb{Z}(p^\infty))$.  We will show that this theory is not categorical as an elementary class, but it is categorical in all cardinals (and more) in the class of abelian torsion groups\footnote{More formally, this means that the class $Mod(Th(\oplus \mathbb{Z}(p^\infty)))$ is not categorical, but $Mod(Th(\oplus \mathbb{Z}(p^\infty)) \cup \{\forall x \bigvee_{n<\omega} n \cdot x = 0\})$ is categorical in all cardinals.}.

This gives a concrete example of a torsion group where more stability theoretic machinery is available when viewing it as a member of a nonelementary class.

For the first part, we note the following general fact.

\begin{prop}
If $R$ is a PID and $M$ is a torsion module such that 
\begin{enumerate}
	\item $ann M = \{0\}$; and
	\item there is $r\in R$ such that $\{m \in M : r \in \mathcal{O}(m) \}$,
\end{enumerate}
then $Th(M)$ is not categorical in any $\lambda \geq |R|$.
\end{prop}

If $R = \mathbb{Z}$, then the first condition says that $M$ is not of finite exponent.

{\bf Proof:} WLOG, $|M| = |R|$.  The given conditions ensure that, for any torsion $M' \equiv M$,
$$M' \precneqq \prod^{tor} M'/U \precneqq \prod M/U$$
and that $\prod M/U$ is not torsion.  Thus, we have a torsion and non-torsion module elementarily equivalent to $M$ in all cardinalities of size at least $|R|$.  Since torsion and non-torsion modules are obviously non-isomorphic, we have the result.\hfill \dag\\

For the second part, we show that, if we have torsion $G \equiv \oplus_{n<\omega} \mathbb{Z}(p^\infty)$, then $G \cong \oplus_{i < |G|} \mathbb{Z}(p^\infty)$.  We rely on the following well-known fact about divisible abelian groups.

\begin{fact}%\cite{?}
Every divisible group is isomorphic to a direct sum of copies of $\mathbb{Q}$ and $\mathbb{Z}(q^\infty)$.
\end{fact}

Since we know $G$ is a divisible $p$-group, it cannot contain any copies of $\mathbb{Q}$ or $\mathbb{Z}(q^\infty)$ for $q \neq p$.  Also, every element of $G$ has infinitely many $p$th roots, so it cannot be a direct sum of finitely many copies of $\mathbb{Z}(p^\infty)$.  Thus, $G \cong \oplus_{i < |G|} \mathbb{Z}(p^\infty)$.

\subsection{Some Stability Theory} \label{stab-sec}

Now that a compactness result is established, we wish to explore some stability theory for torsion modules over a PID (considered as an AEC).  It is already known that all modules are stable (see \cite[Theorem 3.A]{prestmodulebook}), but we have seen that examining nonelementary classes can give stronger results.

\begin{defin}
Let $M$ be an infinite torsion module over a PID.  Then $K_M$ is the AEC whose models are all torsion modules elementarily equivalent to $M$ (in the sense of first-order logic) and where $\prec$ is elementary substructure.
\end{defin}

It is easy to see that this is an AEC.  This class is averageable by Theorem \ref{fulllos-thm}.  Furthermore, this AEC is nicely behaved in the sense that amalgamation holds; Galois types are (first-order) syntactic types; and $K_M$ has no maximal models holds unless all models are isomorphic to $M$.

\begin{prop}
$K_M$ has amalgamation, joint embedding, and Galois types are syntactic.
\end{prop}

{\bf Proof:} We show amalgamation and Galois types are syntactic together.  Note that, since strong substructure is elementary substructure, having the same Galois type implies having the same syntactic type.  Let $M_0 \prec M_1, M_2$ from $K_M$, possibly with $a_\ell \in M_\ell$ such that $tp(a_1/M_0; M_1) = tp(a_2/M_0, M_2)$.  Then, by amalgamation for first order theories, we can find a module $N^*$ and $f_\ell :M_\ell \to N^*$, for $\ell = 1, 2$, that agree on $M_0$ and (if they exist) $f_1(a_1) = f_2(a_2)$.  Set $N$ to be the torsion subgroup of $N^*$.  The torsion radical preserves pure embeddings and picks out a pure subgroup, so we have
$$f_\ell(M_\ell) \subset_{pure} N \subset_{pure} N^*$$
Since they are elementarily equivalent, $Inv(f_\ell(M_\ell), \phi, \psi) = Inv(N^*, \phi, \psi)$ for all p. p. formulas $\phi$ and $\psi$ and pureness implies that $N$ has the same invariants conditions.  Thus, $N$ is elementarily equivalent and the embeddings are elementary.  Furthermore, the $f_\ell$ witness that $a_1$ and $a_2$ have the same Galois type.

The proof of joint embedding is similar.\hfill \dag \\

Note that although two tuples having the same Galois type and having the same syntactic type are the same, the ``Galois types are syntactic'' result above should \emph{not} be taken to mean that any consistent, complete set of formulas is realized as a Galois type; obviously, this is not true a non-torsion element.  However, we do have a local compactness result: any partial type is realizable iff there is a fixed order such that all of its finite subsets are realizable with that fixed order.

\begin{prop}
Let $A \subset N \in K_M$ and $p(\bx)$ be a consistent set of formulas with parameters in $A$.  Then there is an extension of $N$ that realizes $p$ iff there are $r_0, \dots, r_{n-1} \in R$ such that, for every finite $q \subset p$, there is an extension of $N$ that realizes $q \cup \{r_i \cdot x_i = 0 : i < n\}$.
\end{prop}

{\bf Proof:}  This is Theorem \ref{gamcomp} in this context.\hfill \dag\\
%For simplicity, we work with singletons.  Clearly the first condition implies the second.  Now, suppose the second condition that there is $r \in R$, for each finite $q \subset p$, there is some $N_q \succ N$ and $a_q \in N_q$ such that $N_q \vDash r \cdot a_q = 0 \wedge \bigwedge q[a_q]$.  Now, take the torsion ultraproduct $\prod_{q \in P_\omega p}^{tor} N_q / U$.  Since each $N_q$ contains $N$, the normal ``ultrapower'' embedding $h$ takes $N$ into $\prod_{q \in P_\omega p}^{tor} N_q / U$.  Also, since they have fixed order, $a_* := [q \mapsto a_q]_U$ is in this ultraproduct.  The standard argument shows that (the pullback of) $a_*$ realizes $p$.\hfill \dag \\

\begin{prop}
$K_M$ has no maximal models or consists of a single model up to isomorphism.
\end{prop}

{\bf Proof:} Follows directly from Theorem \ref{divline2-thm}. \hfill \dag\\

Since this class is stable, we have a unique independence relation.

\begin{theorem}
\begin{enumerate}
	\item $K_M$ is Galois stable at least in all $\lambda^\omega$.
	\item Coheir is a stability-like independence relation.
	\item Any independence relation satisfying Existence, Extension, and Uniqueness is coheir.
\end{enumerate}
\end{theorem}

{\bf Proof:} \begin{enumerate}
	\item $K_M$ has less types than $Mod(Th(M))$ (this uses that Galois types are syntactic), which is stable (see \cite[Chapter 3, Example 1]{prestmodulebook}).
	\item This follows from Theorem \ref{5.1t}.
	\item This is \cite[Corollary 5.18]{bgkv}
\end{enumerate}

As a consequence of the last statement, this means that the good frame defined by nonsplitting from Vasey \cite{sebframe} is the same as coheir in $K_{\oplus_{n<\omega} \mathbb{Z}(p^\infty)}$.

\end{document}